\documentclass[10pt]{article}
\usepackage{amsmath}
\usepackage{latexsym}
\usepackage{amssymb}
\usepackage{amsfonts}
\usepackage{amsthm}
\title{DIFFERENTIAL GALOIS THEORY AND \\  HOPF ALGEBRAS FOR LIE PSEUDOGROUPS}
\author{J.-F. Pommaret \\ CERMICS, Ecole des Ponts ParisTech, France \\
 jean-francois.pommaret@wanadoo.fr \\
 ORCID: 0000-0003-0907-2601}
\date{  }
\begin{document}
\maketitle

\noindent
{\bf ABSTACT}  \\

According to a quite clever but never acknowledged  work of E. Vessiot that won the prize of the Acad\'{e}mie des Sciences in 1904, " Differential Galois Theory " (DGT) has mainly to do with the study of " Principal Homogeneous Spaces " (PHS) for finite groups ( classical Galois theory), algebraic groups (Picard-Vessiot theory) and algebraic pseudogroups (Drach-Vessiot theory). The corresponding automorphic differential extension are such that $ dim_K(L)< \infty $, transcendence degree $ trd(L/K)< \infty $ and $ trd(L/K)=\infty $ with $ diff trd(L/K)< \infty $ respectively. The purpose of this paper is to mix differential algebra, differential geometry and algebraic geometry in order to revisit DGT, pointing out the deep confusion between {\it prime differential ideals} ( Defined by J.-F. Ritt in 1930) and {\it maximal ideals} thad has been spoiling the works of Vessiot, Drach, Kolchin and all followers. In particular, we use Hopf algebras in order to study the structure of the algebraic Lie pseudogroups involved, namely Lie pseudogroups defined by systems of algebraic OD or PD equations. Many explicit examples are presented for the first time in order to illustrate these results. This paper is also paying a tribute to Prof. A. Bialynicki-Birula on the occasion of his recent death in April 2021 at the age of 90 years old. His main idea has been to notice that an algebraic group $G$ acting on itself is the simplest example of a PHS. If $G$ is defined over a field $K$ and we introduce the algebraic extension $L=K(G)$, then there is a Galois correspondence between the intermediate fields $K \subset K' \subset L$ and the subgroups  $e \subset G' \subset G $  provided that $K'$ is stable under a Lie algebra $\Delta $ of invariant derivations of $L/K$. Our purpose is to extend this result from algebraic groups to algebraic pseudogroups without {\it any} way to use parameters.\\

\vspace{2cm}

Lie pseudogroup; Principal homogeneous space; Rings with operators; Reciprocal distributions.     \\
\vspace{2cm}

\noindent
{\bf COMMENT} \\
The comparison of the paper written by J. Kovacic in 2005 with our second Gordon and Breach book published in 1983 needs no comment on the anteriority of using Hopf algebras in DGT. (See arXiv: 1710.08122, now  published in the book "New Mathematical Methods for Physics", NOVA Science Publishers, New York, 2018).   \\

\newpage

\noindent
{\bf 1) INTRODUCTION}  \\

Before describing the content of this paper, we start providing a few personal historical details explaining the motivation for writing it. We separate this historical background into three parts:  \\

\noindent
A) The story started in $1970$ when the author of this paper was a visiting student of D.C. Spencer in Princeton university. As the mathematical library was opened day and night while being well furnished in french publications, he had the chance to discover the original publication of M. Janet in $1920$ ([15]) on systems of ordinary differential (OD) or partial differential (PD) equations and of E. Vessiot in $1903$ ([45]) on the so-called (at that time) finite and infinite groups of transformations, now respectively called Lie groups and Lie pseudogroups of transformations. This last reference had also been used one year later by Vessiot in order to study the " {\it Differential Galois Theory} " (DGT) in a short but clever paper published in $1904$ ([46]) introducing for the first time " {\it Principal Homogeneous Spaces} " (PHS) for algebraic pseudogroups defined by {\it automorphic systems} of OD or PD equations and the corresponding "{\it automorphic differential extensions} " . The consideration of these results furnished the main content of the first Gordon and Breach (GB) book of $1978$ ([27]). However, he also discovered a rarely quoted long paper published as a Doctor Thesis by J. Drach in $1898$ ([10]), the interesting fact being that the jury was made by G. Darboux, E. Picard and H. Poincar\'{e}, the three most famous french mathematicians existing at that time, a good reason for reading it {\it in the details}. \\

In a rough way, one may say that Vessiot divided his study of DGT into three parts of increasing difficulty, each one being illustrated by original examples. Using standard notations of differential algebra, if $K$ is a differential field of characteristic zero, that is $\mathbb{Q} \subset K$ with $n$ derivations ${\partial}_i$ for $i= 1, ..., n$ and $L$ with derivations $d_i$ for $i=1, ..., n$, in such a way that $K \subset L$ and $d_i {\mid} _K = {\partial}_i$.  \\

\noindent
{\bf a) CLASSICAL GALOIS THEORY}: ${dim}_K L=\mid L/K \mid < \infty$. \\
In this case, $L$ being a finite dimensional vector space over $K$, this means that any element of $L$ is algebraic over $K$ and we refer the reader to the work of E. Galois who died in a duel in $1821$ and to the extensive literature on this subject. The purpose is roughly to study systems of algebraic equations invariant by finite groups, for examples groups of permutations.   \\

\noindent
{\bf b) PICARD-VESSIOT THEORY}: $trd(L/K)< \infty, diff trd(L/K)=0$.  \\
In this case, the maximum number of elements of $L$ which are {\it not} algebraic over $K$ or {\it transcendence degree} of $L$ over $K$, is finite but any element of $L$ must satisfy {\it at least one}  OD or PD algebraic equation over $K$. The problem is roughly to study systems of algebraic OD or PD equations that are invariant by an algebraic group of matrices. The application to the study of shells, chains and analytical mechanics is quite recent ([36]).   \\

\noindent
{\bf c) DRACH-VESSIOT THEORY}: $trd(L/K)=\infty, diff trd (L/K) < \infty$.  \\
In this case, the maximum number of elements of $L$ that are {\it not} algebraic over $K$ may be infinite but one can find a well defined intrinsic maximum number of elements of $L$ that can be given arbitrarily, such a set being called a  {\it differential transcendence basis}. The problem is thus to study systems of algebraic OD or PD equations that are invariant by an algebraic pseudogroup, namely a group of transformations solutions of a system of algebraic OD or PD equations that can be of high order. The best example is provided by the study of the various reductions of the group of invariance of the {\it Hamilton-Jacobi equation} depending on the hamiltonian function of the mechanical system considered, that can be described by subpseudogrous of the pseudogroup of contact transformations, the most important one being the one made by unimodular such transformations ([36],[38]).  \\

The clever idea of Vessiot has been to prove that, {\it in all these cases}, the corresponding Galois theory is just a theory of PHS respectively for finite groups, Lie groups or Lie pseudogroups. Of course, when a Lie group $G$ is acting on a manifold $X$ with an action law $X \times G \rightarrow X: (x,a) \rightarrow y=ax=f(x,a)$, one can introduce the corresponding graph $X \times G \rightarrow X \times X:(x,a) \rightarrow (x,y)$ for the action and $X$ is a PHS for the action of $G$ when such a graph is an isomorphism. The main idea of Vessiot has been to extend such a definition to systems of OD or PD equations and Lie pseudogroups $\Gamma \subset aut(X)$ made by invertible transformations in such a way that, whenever $y=f(x)$ and $\bar{y}=\bar{f}(x)$ are two solutions of this system, then one can find {\it one and only one} invertible transformation $\bar{y}=g(y)$ such that $\bar{f}=g \circ f$ at least locally and such a system is called an {\it automorphic system} by Vessiot. Of course the difficulty, not solved by Drach or Vessiot, was to establish an effective test for checking such a property. In view of the many specific examples presented, the first and second criteria for automorphic systems in Section $3$ will be {\it absolutely unavoidable} in order to decide whether a given differential extension will be automorphic or not (See [28] or [36]) for details).  \\

\noindent
{\bf B)} Having at that time no idea for using the " {\it Spencer machine} " in mathematical physics (See Conclusion !), the author tried to translate the work of Drach and Vessiot in a modern differential geometric language and it is at this precise moment that the true difficulties started. Indeed, though the work of Drach was introducing for the first time crucial concepts like differential fields and other founding stones of differential algebra only sketched by J.F. Ritt more than ... thirty years later ([39]), the author did not succeed, after months of hard work, to understand the thesis of Drach or even to justify the so-called " {\it reductibility} " concept which is the heart of the Drach-Vessiot theory. Even worst, it became clear that the thesis was based on a fundamental confusion between {\it maximal ideal} and {\it prime ideal}, this last concept being only introduced by Ritt in $1930$ ([39]). In order to sketch this difficulty, let us consider the following elementary example that will be revisited later on.  \\

\noindent
{\bf EXAMPLE 1.1}: With $n=1$ independent variable $x$ and $m=1$ dependent variable $y$, let us consider the second order linear system $y_{xx}=0$ while using standard jet notations ([27, 30, 32, 35, 36]). With $K=\mathbb{Q}$ a trivial differential field, we may introduce at once the linear differential ideal $\mathfrak{p} \subset K\{ y \}= K[y, y_x, y_{xx}, ... ] $ which is of course trivially a prime differential ideal generated by the differential polynomial $P=y_{xx}$ and providing successively $y_{xx}=0, y_{xxx}=0, ...$. As a byproduct, taking the full ring of quotients, we may define the differential extension $L= Q(K\{ y \} / \mathfrak{p})= K(y, y_x)$ with $ d_x y = y_x, d_x y_x=0$. It is at once seen to be a differential automorphic extension for the affine group of the real line, namely $ \bar{y} = a^1 y + a^2, {\bar{y}}_x=a^1 y_x$ with two infinitesimal generators $\Theta = \{ {\theta}_1= \frac{\partial}{\partial y}, {\theta}_2 = y \frac{\partial}{\partial y} + y_x \frac{\partial}{\partial y_x}$ and ordinary bracket $[ {\theta}_1 , {\theta}_2 ] = {\theta}_1$. Indeed, we have $a^1 = \frac{{\bar{y}}_x}{y_x} , a^2 = \bar{y} - \frac{{\bar{y}}_x}{y_x} y$ in such a way that the graph of the action $ ((y, y_x), (a^1,a^2)) \rightarrow ((y, y_x), ( \bar{y}, {\bar{y}}_x))$ is an isomorphism. Now, setting $ \Delta = \{ {\delta}_1= y_x \frac{\partial}{ \partial y} , {\delta}_2= y_x \frac{\partial}{\partial y_x} \}$, we obtain $[ \Theta , \Delta ] =0$ that is each $\theta \in \Theta$ commutes with each $\delta \in \Delta$. Accordingly, we may find intermediate fields  $K \subset K' \subset L$ that are {\it not } differential fields because they are not stable under ${\delta}_1$ and are not fields of invariants of subgroups while others are. For example, $K' = \mathbb{Q} ( y y_x)$ is an intermediate field but not a differential field because $ d_x (y y_x)= (y_x)^2\notin K'$. The subgroup of invariance is $\bar{y}= \pm y$ that lets invariant $K'' = \mathbb{Q} (y^2)$ with the strict inclusion $K' \subset K''$. In fact, ${\delta}_1 (y y_x)=(y_x)^2$ and $(y_x)^2/(y y_x)=y_x / y$ on one side, that is $\bar{y} = a^1 y$ but we have also $(y_x)^2 / (\frac{y_x}{y})^2= y^2$ that is $a^1=\pm 1$. On the contrary, $ K' = \mathbb{Q}(\frac{y_x}{y})$ is indeed a differential field because $d_x(\frac{y_x}{y})= - (\frac{y_x}{y})^2 \in K'$ which is also stable by $\Delta$ because ${\delta}_1 (\frac{y_x}{y}) = - (\frac{y_x}{y})^2, {\delta}_2 (\frac{y_x}{y})=\frac{y_x}{y}$ and the biggest group of invariance is $\bar{y}=a y$ that lets invariant {\it exactly} $K'$ in a coherent way with ([3], Theorem 3). We let the reader treat the more general situation with the chain of strict inclusions of differential fields with $\mathbb{Q}<y>= \mathbb{Q}(y,y_x,y_{xx}, ... )$: \\
\[ \fbox{  $  \mathbb{Q} \subset \mathbb{Q} < \frac{y_{xxx}}{y_x} - \frac{3}{2} ( \frac{y_{xx}}{y_x})^2 > \subset \mathbb{Q} < \frac{y_{xx}}{y_x} > \subset \mathbb{Q} < \frac{y_x}{y} > 
\mathbb{Q} < y >   $  } \]  \\

\noindent
{\bf EXAMPLE 1.2}: When $K$ is a differential field containing $\mathbb{Q}$ and $n=1, m=2$, setting $K<y^1,y^2 > = K (y^1,y^2, y^1_x, y^2_x, ... )$ and $K < \Phi >= K (\Phi, d_x \Phi, d_{xx} \Phi, ... )$, for any $\Phi \in K< y >$, we may consider the chain of differential extensions $ \mathbb{Q} \subset \mathbb{Q} < y^2 y^1_x > \subset \mathbb{Q} < y^1, y^2 >$. Using the chain rule for derivatives on the jet level, we obtain easily that the biggest Lie pseudogroup preserving $\Phi = y^2 y^1_x$ is 
$ \Gamma= \{ {\bar{y}}^1 = g(y^1), {\bar{y}}^2 = y^2 / \frac{\partial g}{\partial y^1} \}$ which is defined by the Pfaffian system $ {\bar{y}}^2 {\bar{y}}^1_x= y^2 y^1_x \Rightarrow {\bar{y}}^1 \wedge d {\bar{y}}^2 = d y^1 \wedge d y^2$, namely:   \\
\[   {\bar{y}}^2 \frac{\partial {\bar{y}}^1}{ \partial y^1 }= y^2,  \frac{\partial {\bar{y}}^1}{\partial y^2}=0, \Rightarrow 
\frac{\partial ({\bar{y}}^1, {\bar{y}}^2)}{\partial (y^1, y^2 )} = 1  \Rightarrow \frac{1}{{\bar{y}}^2 }\frac{\partial {\bar{y}}^2}{ \partial y^2}= \frac{1}{y^2}     \]
and we obtain therefore a well defined algebraic pseudogroup. Considering the chain of differential fields $ K =\mathbb{Q} < y^2 y^1_x > \subset K' = \mathbb{Q} < y^2 y^1_x, y^2_x > \subset \mathbb{Q} < y^1, y^2 > = L $, we notice that the biggest pseudogroup preserving $K'$ is the subpseudogroup $ {\Gamma}' = \{ {\bar{y}}^1= y^1 + a, {\bar{y}}^2= y^2 \} \subset \Gamma $ which is preserving the differential field $K'' = \mathbb{Q} < y^2 y^1_x, y^2 > = \mathbb{Q} < y^1_x, y^2 >$ with the strict inclusion $K' \subset K'' $. Once again, it does not seem possible to establish a Galois correspondence between intermediate differential fields and subpseudogroups of $\Gamma$.  \\

\noindent
{\bf C) } It is at this precise moment ($1979$) that the author of this paper discovered that M. Janet was still alive, living a few blocks away in Paris {$16$}° and in a rather good condition as he died quite later on in $1983$ at the age of $96$. He was a close friend of Vessiot and became quite pleased by the fact that the $1978$ GB book had been dedicated to him, for his $88$th Birthday. As a byproduct, he gave him a bunch of private documents now put as a deposit in the main Library of Ecole Normale Sup\'{e}rieure (ENS Paris) where they can be consulted ([29, 38]). We quote below a few lines from a letter from P. Painlev\'{e} to Vessiot  about the thesis of J. Drach (Paris, October 17, 1898):  \\
" Dear Friend, I just read Drach's thesis and I agree entirely with you about the inaccuracy of the two fundamental theorem and their proofs. The mistake is so big that I can hardly conceive that it has been overlooked by the author and the jury. .... I said to Picard that I had received a letter from you about this question; after a few minutes of explanation, he was asthonished to have missed this. I do believe, indeed, there cannot be the slightest shadow of a doubt for anybody making his mind on this thing ... " (See [29], Appendix for the original letters given by Janet). \\ 
Then a " {\it Mathematical Affair} " started in a rather unpleasant way because Drach, supported by E. Borel, never accepted to have made a mistake. Janet also confirmed that the {\it only} paper written by E. Cartan on PD equations, taken from his letters to A. Einstein on {\it Absolute Parallelism} in $1930$ ([4]), never indicated that it was straightly coming from the classical similar paper of Janet  published ... in $1920$ ([15]), without quoting his source of course. Similarly, Cartan {\it never} (we insist !) said that the {\it Vessiot structure equations}  could be competing or even superseding the {\it Cartan structure equations} though both had been created about at the same time around $1905$. {\it Last but not least}, Janet confirmed that everybody in $1930$ knew that the future of DGT should pass through algebraic pseudogroups and {\it not} (I insist again) through the so-called " {\it differential groups} " of J. F. Ritt ([5, 18, 26, 39]) that were considered as leading to a dead end because only trivial examples had been exhibited.  \\

As a byproduct, after loosing almost half a year for nothing, the author lectured during a month at the Columbia University  in New York along an invitation of E. Kolchin who, even after discovering these facts, refused to provide a preface for the DGT book of $1983$ ([28]). It is easy for any reader to check along the list of the differential algebra community provided by J. Kovacic, including himself as a former student of Kolchin at that time, that not a single of these more than hundred persons did even quote {\it once} this DGT book since ... about forty years, despite the fact that Hopf  algebras had been used for the first time in this book ([19, 20]). Finally, to be fully fair, we may say that, being invited to lecture at the King's College of London about the Backlund problem and jet theory by F. A. E. Pirani, D. C. Robinson and W.F. Shadwick just after the publication of the $1978$ GB book, it happened that the only visitor's room left was the one occupied by M. E. Sweedler ([6, 44]), on leave for a week, where bunches of papers on Hopf algebras were on display on the table, a huge chance indeed. \\

We finally explain on a simple example the content of the paper written by A. Bialynicki-Birula in $1961$ ([2]) and using for the first time tensor products of rings and fields for studying DGT.  \\
 {\it The main idea is that a Lie group acting on itself is the simplest example of a PHS}. \\

\noindent
{\bf EXAMPLE 1.3}: Let us consider again the group of affine transformations on the real line. If $y=a^1 x + a^2, z=b^1 y + b^2$, we obtain by composition $z= (b^1 a^1) x + (b^1 a^2 + b^2)$ and the group $G$ with composition $({\bar{a}}^1, {\bar{a}}^2)= (b^1,b^2) (a^1, a^2)=( b^1 a^1 , b^1 a^2 + b^2)$ and inversion $(a^1, a^2 )^{-1}= \frac{1}{a^1}, - \frac{a^2}{a^1}) $. The two reciprocal commuting left and right distributions on $G$ are respectively generated by: \\
\[  \fbox{  $   \begin{array}{lcl}
\Theta & = & \{   {\theta}_1= a^1 \frac{\partial}{\partial a^1} + a^2 \frac{\partial}{\partial a^2} , {\theta}_2= \frac{\partial}{\partial a^2}\} \\
            &    &                                \\
\Delta & = & \{ {\delta}_1 = a^1 \frac{\partial}{\partial a^1} , {\delta}_2 = a^1 \frac{\partial}{\partial a^2}  \} 
\end{array}    $  }  \]
in such a way that $[ \Theta , \Delta ]=0 $, that is $ [ {\theta}_i , {\delta} _j ] =0, \forall i, j =1, 2 $. 
Of course, the group $G$ being defined on $k = \mathbb{Q}$, we have $ k[G]=k [a^1, a^2 ]$ and we may introduce the fields 
$ K = \mathbb{Q}  $ , $L = \mathbb{Q} (a^1, a^2) $ with $K \subset L $.               
Introducing the intermediate field $ K' = \mathbb{Q}(a^2/a^1)$ with $K \subset K' \subset L$, we obtain at once:  \\
\[      \frac{b^1 a^2 + b^2}{b^1 a^1} = \frac{a^2}{a^1} + \frac{b^2}{b^1 a^1} = \frac{a^2}{a^1}   \Longleftrightarrow   b^2=0  \]
The subgroup preserving $K'$ has parameters $(b^1,0)$ in such a way that $(b^1,0) (a^1, a^2)= b^1 a^1, b^1 a^2) $ and the only invariant subfield is again $K'$. On the contrary, choosing $K'=\mathbb{Q}(a^1 a^2)$, we obtain:  \\
\[      (b^1 a^1)(b^1 a^2 + b^2) = (b^1)^2 (a^1 a^2) + (b^1 b^2 a^1) = a^1 a^2 \longleftrightarrow (b^1)^2=1, b^2 = 0 \]
and the invariant subfield is now $K'' = \mathbb{Q}((a^1)^2, a^1 a^ 2, (a^2)^2)$ with the strict inclusion $ K \subset K' \subset K'' $. \\
BB did notice that only the intermediate fields stabilized by $\Delta$ do provide a Galois correspondence $K \subset K' \subset L \leftrightarrow e \subset G' \subset G$ as we have indeed ${\delta}_1 (\frac{a^2}{a^1}) = 0,\,\,\,  {\delta}_2 (\frac{a^2}{a^1}) = 1 $ but:
\[ {\delta}_1 (a^1 a^2) = a^1 a^2 , \,\,\, {\delta}_2 (a^1 a^2) = (a^1)^2  \Rightarrow (a^2)^2 = (a^1 a^2 )^2 / (a^1)^2       \]

To recapitulate, in his first paper BB discovered that only intermediate fields stable under $\Delta$ do provide a Galois correspondence in general. Such a comment is at the origin of his second paper providing a new approach to the Picard-Vessiot theory by using "{\it fields with operators}", differential fields being only a specific example. In a rough way, his first idea has been to introduce the tensor product $L {\otimes}_K L$ with $L=\mathbb{Q}(a^1,a^2) $ on the left and $ L= \mathbb{Q}({\bar{a}}^1, {\bar{a}}^2)$ on the right. His second idea has been to enlarge $\Delta$ to the new derivations {\it on the tensor product}, namely:  \\
\[ \fbox{  $  \Delta = \{ {\delta}_1 = a^1 \frac{\partial}{\partial a^1} + {\bar{a}}^1 \frac{\partial}{\partial {\bar{a}}^1}, \,\,
                     {\delta}_2 = a^1 \frac{\partial}{\partial a^2} + {\bar{a}}^1 \frac{\partial}{\partial {\bar{a}}^2}  \}  $   }   \]
He then noticed that the only quantities killed by the new ${\delta}_1$ and ${\delta}_2$ are:  \\
\[  \fbox{ $   \{  \,\,\,  b^1 = \frac{{\bar{a}}^1}{a^1} = (\frac{1}{a^1}) \otimes {\bar{a}}^1, \hspace{1cm}   
b^2 = {\bar{a}}^2 - (\frac{{\bar{a}}^1}{a^1} ) a^2 = 1 \otimes{\bar{a}}^2 - (\frac{a^2}{a^1}) \otimes {\bar{a}}^1 \,\,\,   \}  $  } \]
It follows that $\mathbb{Q}[b^1,b^2] =\mathbb{Q} [G] = cst (L {\otimes}_K L) $ where the “ {\it constants} " are now the quantities killed by the extension of the derivations $\Delta$ to $L {\otimes}_K L$. We do consider that the upper " {\it bar} " notation is superseding the tensorial notation in the differential geometric framework. The novelty (we could even say revolution) is that the group parameters are no longer " {\it differential constants} " contrary to the standard approach of Kolchin that can be sketched as follows on Example 1.1 as we get indeed $(a^1 = \frac{{\bar{y}}_x}{y_x}, a^2= \bar{y} - y\frac{{\bar{y}}_x}{y_x}) \Rightarrow  (d_x a^1 = 0, d_x a^2 = 0)$ because $y_{xx}=0 \Rightarrow {\bar{y}}_{xx}=0$.  \\

The content of the paper is now clear from this long Introduction. In the second section, we revisit the classical Galois theory. In the third section, we recall the definition of rings and corings in a purely algebraic framework. In the fourth section, we present the basic concepts of jet theory and Lie pseudogroups in order to understand the concepts of automorphic system and automorphic differential extension. Finally, in the fifth section, we provide motivating examples in order to illustrate the main results obtained by using Hopf algebras for the differential Galois theory (DGT), before the conclusion in the sixth section. \\

We hope that the comparison of our two books ([28, 36]) with the paper published by J. Kovacic ($1941$-$2009$) in $2005$ ([20]) needs no comment on the anteriority of using Hopf algebras in DGT. We also point out the fact that it is difficult to understand the conceptual confusion done during almost fifty years between algebraic pseudogroups and differential algebraic groups, knowing that such a comment had already been made in $1930$. We believe that the reason is mainly coming from the fact that Ritt or Kolchin were involved in analysis rather than in differential geometry, always caring about solutions of systems of OD or PD equations. In any case, future will judge !.  \\  \\  \\

\noindent
{\bf 2) CLASSICAL GALOIS THEORY REVISITED}. \\

Let us start this section explaining the {\it clever idea} of E. Vessiot in the first chapter of his $1904$ paper ([46]), both with the {\it deep confusion} that is spoiling it. The following example will prove that the standard link existing between the classical Galois theory and group theory is {\it not at all} as evident as one could imagine from the extensive literature on the subject ([1, 7, 11, 23, 24, 12, 41, 43]) and the fact that tensor products of rings and fields have rarely been used ([14, 22, 47]). \\

\noindent
{\bf EXAMPLE 2.1}: First of all, with ground field $K=\mathbb{Q}$ and one indeterminate $y$, in order to understand the distinction that may exist between {\it general} and {\it special} algebraic equations, let us consider the {\it general} cubic polynomial equation $P \equiv y^3 -{\omega}^1 y^2 + {\omega}^2 y - {\omega}^3 = 0$ with upper indices on $\omega$ in order to agree with the next sections. In algebra, giving {\it special} values in $K$ to $\omega$, the main problem has always been to know about the three roots $({\eta}_1, {\eta}_2, {\eta}_3)$, that is to say to exhibit a {\it splitting field} $K({\eta}_1, {\eta}_2, {\eta}_3)$ of $P$. Comparing now $P$ to the product $ (y-{\eta}_1)(y-{\eta}_2)(y-{\eta}_3)$, we obtain the three symmetric functions of the roots: 
\[ \fbox{ $   {\eta}_1 + {\eta}_2 + {\eta}_3 = {\omega}^1, \,\,\, {\eta}_1{\eta}_2 + {\eta}_1{\eta}_3 + {\eta}_2 {\eta}_3 = {\omega}^2, \,\,\,  {\eta}_1 {\eta}_2 {\eta}_3 = {\omega}^3  $ } \]
The roots are different if and only if $\delta = ({\eta}_1 - {\eta}_2)({\eta}_1 - {\eta}_3)({\eta}_2 - {\eta}_3) \neq 0$. Introducing the derivative $P'=dP/dy$, it is easy to construct the {\it resultant} of $P$ and $P'$ which only depends on $\omega$ and is equal to ${\delta}^2$ up to a factor in $\mathbb{Q}$. Modifying slightly standard notations of textbooks, we shall obtain:
\[  \fbox{ $  {\delta}^2= - 27 ({\omega}^3)^2 + 18 {\omega}^1 {\omega}^2 {\omega}^3 - 4 ({\omega}^2)^3 - 4 ({\omega}^1)^3{\omega}^3 + ({\omega}^1{\omega}^2)^2 $  }  \]
We invite the reader to consider the special polynomial $P \equiv y^3 - 3 y + 1 = 0$ before reading ahead in order to discover that it is not so natural to associate as Galois group the permutation subgroup:
 \[A_3 = \{ e=(\substack{123\\123}), (\substack{123\\231}), (\substack{123\\312}) \} \lhd S_3= \{ e=(\substack{123\\123}), (\substack{123\\231}), (\substack{123\\312}), (\substack{123\\132}), (\substack{123\\321}), (\substack{123\\213}) \} \]
with $\mid A_3 \mid = 3 < \mid S_3 \mid = 6$ when studying the Galois extension $L=Q(K[y]/(P))/K$ generated by the polynomial $P$ which is irreducible over $K$ or, equivalently, $P$ has no root in $K$.  \\

Indeed, let us introduce the generic zero $\eta$ by the specialization $K[y] \rightarrow K[y]/ (P): y \rightarrow \eta$ in order to get ${\eta}^3 - 3 \eta +1 = 0$. Any automorphism $\sigma \in aut(L/K)$ is of the form $\sigma(\eta) : \alpha 1 + \beta \eta + \gamma {\eta}^2$ with $\alpha, \beta, \gamma \in k$. A straightforward but quite tedious substitution left to the reader proves that, if $\eta$ is a root of $P$, then $\sigma (\eta) = {\eta}^2 - 2$ is another root, also with ${\sigma}^2(\eta)= {\eta}^4 - 4{\eta}^2 + 2 = - {\eta}^2 + \eta +2$ and that ${\sigma}^3(\eta)= \eta \Rightarrow {\sigma}^3=e$ is the identity. It follows that the splitting field of $L/K$ is just $L/K$ which is therefore a Galois extension as we already said. The main problem is that the definitions of $\eta$ and $\sigma$ only depend on $\L/K$ and do not provide an {\it algebraic group} "per se". The following "trick" allows to get rid of such a problem ([41]). For this, we notice that the equation can be written as $y^3=3y - 1 \Leftrightarrow y(y^2-2)= y - 1$. We may therefore consider the purely rational transformation $y \rightarrow \sigma(y)=\bar{y}= 1 - \frac{1}{y}$ defined over $\mathbb{Q}$, obtaining thus ${\sigma}^2(y) =  1 - (1 / (1 - \frac{1}{y})) = \frac{1}{1-y}$ and ${\sigma}^3(y)=y$ in a way totally independent of the equation. Such a new way will permit to present the classical Galois theory along the way initiated by Vessiot, through the so-called {\it general equations}. \\

We have a finite algebraic group $\Gamma$ of transformations of the real line $Y$ which is defined over the ground field $k=\mathbb{Q}$. Considering the rational function:\\
\[    y \rightarrow  \Phi (y)\equiv  y + \sigma (y) + {\sigma}^2(y) = y + (1 - \frac{1}{y}  ) + ( \frac{1}{1-y}) = \frac{y^3 - 3 y + 1}{y^2 - y} \in k(y) = L  \]
Let us prove that such a function, which is of course invariant by $\sigma$, is the generating invariant of the rational action of $\Gamma$ on $k(y)$ as follows. Writing $\Phi(\bar{y})=\Phi(y)$, we get successively: 
\[    (y^2 - y)({\bar{y}}^3 - 3 \bar{y} + 1)- ({\bar{y}}^2 - \bar{y})( y^3 - 3 y + 1)=0  \,\, \Rightarrow \,\,     (\bar{y} - y)(\bar{y} - (1 - \frac{1}{y}))(\bar{y} - \frac{1}{1-y})  =0  \]
With $k=\mathbb{Q}, K=k(\Phi), L=k(y)$, we obtain the chain of inclusions $k \subset K \subset L$ with $\mid L/K \mid =3$ and $k$ algebraically closed in $L$, thus in$K$. We discover that $Y$ is a {\it Principal Homogeneous Space} (PHS) for $\Gamma$ with graph $Y \times Y \simeq Y\times \Gamma$ with a sight abuse of language but with a well defined isomorphism $L{\otimes}_K L \simeq L {\otimes}_k k[\Gamma] $ which can be extended to the rings of quotients $Q(L{\otimes}_K L) \simeq Q(L {\otimes}_k k[\Gamma])$ though each side is already a direct sum of $3$ fields as we saw. Introducing what we shall call the {\it Lie form} $\Phi(y)=\omega$, of the action, we finally notice that the left term is defined by using the so-called {\it general} equations $ y^3 - \omega y^2 + (\omega - 3) y +1 = 0$ and ${\bar{y}}^3 - \omega {\bar{y}}^2 + (\omega - 3) \bar{y} +1 =0$. Indeed, subtracting the first 
from the second, we get: 
\[  (\bar{y} - y) (( {\bar{y}}^2 + y \bar{y} + y^2)  - \omega (\bar{y} + y) + \omega - 3) = (\bar{y} - y) ( {\bar{y}}^2 - (\omega - y) \bar{y} - \frac{1}{y})=0  \]
because, dividing by $y$ the first general equation, we get $y^2 - \omega y + \omega - 3 = - \frac{1}{y}$ and obtain the same factorization as before as $(1- \frac{1}{y})+ \frac{1}{1-y}= \omega - y, (1- \frac{1}{y})(\frac{1}{1-y}) = - \frac{1}{y}$.  \\ 
We finally explain on this example how Vessiot got in mind the fact that classical Galois theory is only a theory of PHS for groups of permutations. For this, let us introduce  $3$ indeterminates $y=(y^1, y^2, y^3)$ or $y^k$ for $k=1,2,3$ and consider the $3$ equations $P(y^k)=0$ as a linear system for $({\omega}^1, {\omega}^2, {\omega}^3)$ with a Van der Mond determinant which is homogeneous of degree $3$ in $\mathbb{Q}[y]$, thus equal up to sign to $\Delta(y)=(y^1 - y^2)(y^1 - y^3)( y^2 - y^3)$. We obtain $3$ general equations in Lie form:
\[  {\Phi}^1(y) \equiv y^1 + y^2 + y^3={\omega}^1, \,\, {\Phi}^2 (y)\equiv y^1 y^2 + y^1 y^3 + y^2 y^3 = {\omega}^2, \,\, {\Phi}^3(y) \equiv y^1 y^2 y^3 = {\omega}^3 \]
We have proved in ([28], p 151) that the ideal $\mathfrak{a}$ generated by the equations $\Phi - \omega=0$ is perfect if and only if $\delta\neq 0$, that is when the $3$ roots of the equation $P(y)=0$ are different.   \\
The main definition of Vessiot has been the following:  \\

\noindent
{\bf DEFINITION 2.2 }: "A system of equations is called an {\it automorphic system} if any solution may be obtained from a given one by one and only one transformation of a (eventually finite) group of transformations acting on the variables".  \\

In the present situation, we have indeed a PHS for the group $S_3$ of permutations in $3$ variables. However, if $\mathfrak{a} \subset \mathbb{Q}[y]$ is prime for the general situation and for certain special situations like $for P\equiv y^3 +y +1=0$ with ${\delta}^2= - 31$, it may not be prime for others, in particular the present one $P\equiv y^3 - 3 y +1=0$ with ${\delta}^2 =81=9^2$, because $\mathfrak{a}$ surely contains the product $(\Delta(y) - 9)(\Delta(y) + 9)$ though each of the factors does not belong to $\mathfrak{a}$. Adding the equation $\Delta(y)  = 9$, we get a prime ideal reducing the group of invariance from $S_3$ to $A_3$. However, Vessiot was writing in $1904$ and the concept of a {\it prime ideal} has only been introduced in $1930$ by J.F. Ritt when he created differential algebra ([28, 39]). In a rough way, Vessiot has been confusing "{\it prime ideal}" with "{\it maximal ideal}" (thus prime) in his (personal) definition of {\it irreducibility}. The worst fact is that the whole Picard-Vessiot theory has also been based on this confusion, not known or even acknowledged by E.R. Kolchin, despite what we told him in front of his students while lecturing at the Columbia University of New York ($8$ lectures of $2$ hours in April $1981$! ) ([16, 17]). The reason is also that Kolchin was engaged in a kind of "dead end" with his so-called {\it differential algebraic groups} along with a confusing definition provided by J.-F. Ritt in the last few papers he wrote around $1950$, just before he died ([5, 18, 26]). We invite the reader to treat similarly the general quadratic equation $y^2 - \omega y + 1= 0$ with group $\Gamma = \{ \bar{y}=y, \bar{y}= \frac{1}{y} \}$ while setting $\Phi(y)\equiv y + \frac{1}{y}=\omega$.   \\

In characteristic $0$, let us start recalling a few technical results on tensor products of rings and fields that are not so well know, following closely ([28, 47]). In this section, we shall only deal with finitely generated field extensions, contrary to the next sections.  \\

\noindent
{\bf DEFINITION 2.3}: Let $A$ is a ring with unit $1$ and elements $a,b,c,..$. If $\mathfrak{a}\subset A$ is an ideal, we may introduce its {\it radical} to be $rad(\mathfrak{a})=\{a\in A \mid \exists n \in \mathbb{N}, a^n \in  \mathfrak{a} \}$ which is also sometimes simply denoted by $\sqrt{\mathfrak{a}}$. An ideal $\mathfrak{a}$ is said to be {\it perfect} if $rad(\mathfrak{a})=\mathfrak{a}$ and the residue ring $k[y]/\mathfrak{a}=A$ is said to be {\it reduced} in this case. An ideal $\mathfrak{p}\in A$ is said to be {\it prime} if $ab \in \mathfrak{p} \Rightarrow  a\in \mathfrak{p} \,\, or \,\, b\in \mathfrak{p}$ and the residue ring $A / \mathfrak{p}$ is thus an integral domain in this case. because $\bar{ab} = 0 \Rightarrow \bar{a}=0$ or $\bar{b}=0$ y denoting a residue with a bar.  \\

\noindent
{\bf PROPOSITION 2.4}: Any perfect ideal $\mathfrak{a}$ in a polynomial ring is the intersection $\mathfrak{a}={\mathfrak{p}}_1 \cap ... \cap {\mathfrak{p}}_r$ of a finite number $r$ of prime polynomial idea.4ls. \\

\noindent
{\bf PROPOSITION 2.5}: A maximal ideal $\mathfrak{m} \subset max(A)$ is prime and an ideal $\mathfrak{m}$ is maximal if and only if the residue ring $A/ \mathfrak{m}$ is a field.\\ 

\noindent
{\it Proof}: If $ab \in \mathfrak{m}, a \notin \mathfrak{m}$ then $\mathfrak{m} + a A =A \Rightarrow \exists c \in A, ac - 1 \in \mathfrak{m} \Rightarrow b \in \mathfrak{m}$. Then, let us first prove that $\mathfrak{m}\in max(A) \Rightarrow  A / \mathfrak{m}$ is a field. If $a \in A, a\notin \mathfrak{m}$, then $\exists$ an ideal $\mathfrak{a}= \{ a b + c \in A \mid b \in A, c \in \mathfrak{m} \} \subset A $. Indeed, we have successively $u(ab+c)=a(ub) + uc \in \mathfrak{a}, (ab_1 + c_1) + (ab_2 + c_2) =a(b_1 + b_2) + (c_1 + c_2) \in \mathfrak{a}$. Also, choosing $b=0, c\in \mathfrak{m}$, we get $\mathfrak{m} \subset \mathfrak{a}$ and a contradiction unless $\mathfrak{a} = A$. We may thus find $b$ and $c$ such that $ab+c=1$. Passing to the residue $A \rightarrow A/ \mathfrak{m}: a \rightarrow \bar{a}$, we get $\bar{a}\neq 0$ and $ \bar{a} \bar{b}=1$, that is $A/ \mathfrak{m}$ is a field. \\
{\it Conversely}, let us imagine that $A / \mathfrak{p}$ is a field, that is $\mathfrak{p}$ is at least a prime ideal, but $\mathfrak{p} \notin max(A)$. Then, we may find ideals $ \mathfrak{p} \subsetneqq \mathfrak{q} \subset A$ and choose $a \in \mathfrak{q}, a\notin \mathfrak{p}$ in such a way that $\mathfrak{p} \subsetneqq a A + \mathfrak{p} \subseteq \mathfrak{q}  \subset A$. Now, as $A / \mathfrak{p}$ is a field, $\exists b\in A, \bar{a} \bar{b} = 1 \Rightarrow (a + \mathfrak{p})(b + \mathfrak{p})= 1 + \mathfrak{p} \Rightarrow  \exists c \in \mathfrak{p} , ab + c = 1$ and thus $a A +  \mathfrak{p} = A \Rightarrow \mathfrak{q}= A \Rightarrow \mathfrak{p} = \mathfrak{m} \in max(A)$. \\
$\hspace*{12cm} \Box$ \\

\noindent
{\bf DEFINITION 2.6}: When $R$ is a ring, the subrings $A$ and $B$ containing a subfield $k$ are said to be {\it linearly disjoint} over $k$ in $R$ if, whenever $a_1, ..., a_r \in A$ are linearly independent over $k$ and $b_1, ...,b_s \in B$ are linearly independent over $k$, then the $rs$ products $a_ib_j$ are linearly independent over $k$ in $R$. We shall denote by $[A,B]$ the smallest subring of $R$ containing both $A$ and $B$. \\
Similarly, when $N$ is a field, the subfields $K$ and $L$ both containing a subfield $k$ are said to be {\it linearly disjoint} over $k$ in $N$ if, whenever $a_1, ..., a_r \in K$ are linearly independent over $k$ and $b_1, ...,b_s \in L$ are linearly independent over $k$, then the $rs$ products $a_ib_j$ are linearly independent over $k$ in $N$. We shall denote by $(K,L)$ the smallest subfield of $N$ containing both $K$ and $L$.\\

The two following propositions will be quite useful along with the following diagram ([2, 3]):  \\

\[ \begin{array}{cccccc}
M &   \longrightarrow  & (K,M) & \longrightarrow   &  (L,M)   \\
\uparrow  &       &   \uparrow  & & \uparrow     \\
k & \longrightarrow     &   K  & \longrightarrow   &  L
\end{array}      \]    \\

\noindent
{\bf PROPOSITION 2.7}: If $k \subset K \subset L$ and $k \subset M$ are subfields of a bigger field $N$, then $L$ and $M$ are linearly disjoint over $k$ in $(L,M)$ if and only if $L$ and $(K,M)$ are linearly disjoint over $K$ in $(L,M)$ {\it and} $K$ and $M$ are linearly disjoint over $k$ in $(K,M)$.  \\

\noindent
{\it Proof}: As vector spaces, let $\{{\lambda}_r\}$ be a basis of $K$ over $k$, $\{ {\mu}_s\}$ be a basis of $L$ over $K$ and $\{ {\nu}_t\}$ be a basis of $M$ over $k$. Then $\{ {\lambda}_r{\mu}_s \}$ is a basis of $L$ over $k$. If $L$ and $M$ were {\it not} linearly disjoint over $k$ in $(L,M)$, we may find linear relations of the form ${\Sigma}_{r,s}
( {\Sigma}_{t} c_{rst} {\nu}_t ) {\lambda}_r {\mu}_s = 0$ that we could write as ${\Sigma}_s ({\Sigma}_{r,t} c_{rst} {\lambda}_r {\mu}_t) {\mu}_s =0$, contradicting the linear disjointness of $L$ and $(K,M)$ over $K$. The converse is similar.  \\
$\hspace*{12cm}  \Box $  \\

\noindent
{\bf PROPOSITION  2.8}: In the situation of the last proposition, we have $(K,M) \cap L = K$.  \\

\noindent
{\it Proof}: Let again $\{ {\nu}_t\}$ be a basis of $M/ k$. Then, any element $\mu \in (K,M)$ may be written as $\mu= {\Sigma}_t a_t {\nu}_t / {\Sigma}_t b_t {\nu}_t$, with $a_t, b_t \in K$. If 
$\mu \in L$, then ${\Sigma}_t (a_t - \mu b_t) {\nu}_t=0$ with $(a_t, \mu b_t) \in L$. However, according to the preceding proposition, $(K,M)$ and $L$ are linearly disjoint over K in 
$(L,M)$ and we must have therefore $\mu b_t = a_t$. Now, one of the $b_t$ {\it at least} must be different from zero because otherwise ${\Sigma}_t b_t {\nu}_t$ could not be used as a denominator. Hence, for some $t$, we have $\mu = a_t / b_t \in K   $. Finally, we have of course $K \subset (K,M), K \subset L \Rightarrow  K \subseteq (K,M) \cap L$ and this ends the proof.  \\
$ \hspace*{12cm}    \Box $  \\

\noindent 
{\bf DEFINITION 2.9 }: An extension $L/K$ is called {\it regular} if $K$ is algebraically closed in $L$.  \\

\noindent
{\bf THEOREM 2.10 }: If $K$ and $L$ are two fields containing a field $k$ and $L/k$ is regular, then $K{\otimes}_k L$ is an integral domain.  \\

\noindent
{\it Proof}: Let us decompose the extension $K/k$ by introducing a transcendence basis $\{s_i\mid i\in I\}$ of $K/k$ in such a way that $K$ becomes algebraic over $k(s) = k(s_1, ..., s_m)$. We have the commutative diagram of inclusions:  \\
\[ \begin{array}{ccccc}
 L  &  \rightarrow    & L(s) & \rightarrow & Q(K{\otimes}_k L)  \\
  \uparrow  & & \uparrow  &   &  \uparrow  \\
k & \rightarrow & k(s) & \rightarrow & K   
  \end{array}      \]
  
By induction on $m$, one can prove ([28], Lemma 4.47) that $L(y)/k(y)$ is regular for any indeterminate $y$ and thus $L(s)/k(s)$ is regular. Also, when $P\in K[y] $ is irreducible over $K$, then it is also irreducible over $L$ when $L/K$ is regular ([28], Proposition 4.48). Hence, we have just to prove that $K {\otimes}_{k(s)} L(s)$ is a field. But this fact follows because $K$ may be generated by a single {it primitive element} element $\eta$ over $k(s)$, generic root of an irreducible polynomial $P\in k(s)[y]$ that remains irreducible over $L(s)$. If we require that $K$ and $L$ be linearly disjoint over $k$ in $(K,L)$, we need that the homomorphism $K{\otimes}_k L \rightarrow (K,L)$ be a monomorphism. In this case there exists an isomorphism $K{\otimes}_k L \simeq [K,L]$, a reason for which we need that $K{\otimes}_k L$ be an integral domain in order to be able to set $[K,L] = K{\otimes}_k L \Rightarrow (K,L) = Q(K{\otimes}_k L)$ {\it exactly}, that is independently of any bigger field $N$ as before. \\
$ \hspace*{12cm}   \Box $ \\

\noindent
{\bf EXAMPLE 2.11 }: Let us consider the automorphic system ([28], Example 8.57, p 177, [    ]):  \\
\[   {\cal{A}}  \hspace{6cm}  \Phi \equiv \frac{y^4 - 1}{y^2} = \omega   \hspace{6cm}  \]
or, equivalently, the general quartic equation $y^4 - \omega y^2 - 1 =0$ which is known to have the Galois group $D_8$ with $\mid D_8 \mid = 8$. For this, considering the new equation ${\bar{y}}^4 - \omega {\bar{y}}^2 - 1 = 0$ and substracting the previous equation, we get $(\bar{y} - y)(\bar{y} + y) ({\bar{y}}^2 + y^2 - \omega)=0$ but the last term is easily seen to be 
$ ({\bar{y}}^2 + \frac{1}{y^2})= (\bar{y} - \frac{i}{y})(\bar{y} + \frac{i}{y})$ and we obtain the finite group $ \Gamma (L/K)=\{ y \rightarrow y, - y , \frac{i}{y}, - \frac{i}{y} \}$ with $\mid L/K \mid = 4$ along with the following diagram in which $split(L/K) = L(i)$ is a smallest Galois extension of $K$ containing $L$ with $\mid split(L/K)/K \mid= \mid D_8 \mid = 8$ in the following picture:
\[     \begin{array}{lcl}
L = \mathbb{Q}(y)& \stackrel{2}{\longrightarrow}  &  L(i) = split(L/K)  \\
\hspace{2mm} \scriptstyle{4}\uparrow   & \scriptstyle {8} \nearrow  & \hspace{3mm} \uparrow \scriptstyle{4}  \\
 K=\mathbb{Q}(\Phi) & \stackrel{2}{\longrightarrow} & K(i)  \\
  \hspace{4mm} \uparrow &   &  \hspace{3mm} \uparrow     \\
k= \mathbb{Q}  & \stackrel{2}{\longrightarrow }  & k(i)
\end{array} \]
Now, $L(i)$ is a Galois extension of $K(i)$ with $\mid L(i)/K(i) \mid = 4$ and we may decompose it into two quadratic extensions. First, we have the subgroup $\{ y \rightarrow y, \frac{i}{y} \}$ with generating invariant $\Psi = y + \frac{i}{y}$, satisfying ${\Psi}^2 - \Phi - 2i = 0$. If we consider the intermediate field $K(i, \Psi)$ between $K(i)$ and $L(i)$, we get successively $K(i, \Psi) \cap L = K \rightarrow (K,k(i))=K(i) \subset K(i,\Psi)$ with a strict inclusion. However, if we consider the other subgroup $\{ y \rightarrow y,  -y \}$, the only generating invariant is $\Psi= y^2$ and we have now $K(i, \Psi) \cap L= K(\Psi) = K' $ . We invite the reader to revisit this example by using a linear group $\Gamma$ of $2 \times 2$-matrices preserving the two invariants  ${\Phi}^1 \equiv y^2 - z^2, {\Phi}^2 \equiv yz$, proving first that:\\
\[ (ay+bz)^2 - (cy+dz)^2 = y^2 - z^2, \,\, (ay+bz)(cy+dz)=yz  \]
provides an ideal $\mathfrak{a} = ( a^2 - c^2 = 1, b^2 - d^2 = - 1, ab - cd = 0, ac = 0, bd = 0, ad + bc = 1 )$ showing that the algebraic group  can be decomposed into the three isolated prime components:
\[  (a= 1, b=0, c=0, d=1) \cup (a= - 1, b=0, c= 0, d= - 1) \cup (a= 0, b^2= - 1, b+c=0, d=0)  \] 
In this case, we have thus $k[\Gamma] \simeq \mathbb{Q} \oplus \mathbb{Q} \oplus \mathbb{Q}(i)$. This example proves that one cannot hope to refer to a classical Galois extension in order to study intermediate fields.  \\
We invite the reader to treat similarly the case ${\Phi}^1 \equiv y^2 + z^2, {\Phi}^2 \equiv yz$ and compare. Equivalently, it amounts to consider the irreducible general equation $y^4 - \omega y^2 +1 =0$ and the Galois extension $L/K$ with $K = \mathbb{Q}(y^2 + \frac{1}{y^2})= \mathbb{Q}( \Phi) \subset L = \mathbb{Q}(y)$ and $\mid L/K \mid = 4$. Then, $L {\otimes}_K L $ s defined by the factorization $(\bar{y} - y)(\bar{y} + y)( {\bar{y}}^2 + y^2 - \omega)$ with last term equal to $({\bar{y}}^2 - \frac{1}{y^2})= (\bar{y}  - \frac{1}{y})(\bar{y} + \frac{1}{y})$. Hence, $\L/K$ is a Galois extension with Galois group the Klein group $V_4$ such that $\mid V_4 \mid = 4$ and thus $k(V_4)\simeq \mathbb{Q}\oplus \mathbb{Q} \oplus \mathbb{Q} \oplus \mathbb{Q}$. The subgroup ${\Gamma}' = \{ \bar{y}=y, \bar{y}= \frac{1}{y} \}$ is defined by the additional invariant $\Psi = y + \frac{1}{y}$ with ${\Psi}^2 - \Phi - 2 =0$. while the other subgroup ${\Gamma}' = \{ \bar{y} = y, \bar{y} = - y \}$ is defined by the additional invariant $\Psi = y^2$ with ${\Psi}^2 - \Phi \Psi + 1 = 0$. In both cases, we have ${\Gamma}' \lhd \Gamma$ and $K'/K$ is a Galois extension.  \\

\noindent
{\bf EXAMPLE  2.12 } : Coming back to the Vessiot point of view, we may choose $k= \mathbb{Q} \subset  K= k({\Phi}^1, {\Phi}^2, {\Phi}^3))=k(\Phi) \subset K'= k(\Phi, \Delta) \subset L= k(y^1, y^2, y^3)= k(y)$ with:  
\[ {\Phi}^1\equiv y^1 + y^2 + y^3, {\Phi}^2 \equiv y^1 y^2 + y^1 y^3 + y^2  y^3, {\Phi}^3 \equiv y^1 y^2 y^3, \Delta \equiv (y^1 - y^2)(y^1 - y^3)(y^2 - y^3) \] 
to revisit the classical Galois theory for $A_3 \lhd S_3$ while understanding that the identification of $\Gamma (L/K)$ with $aut (L/K) $ is just a pure coincidence because Galois extensions are just examples of automorphic extensions for groups of permutations represented by square matrices with entries in $(0, 1)$. Of course, it also remains to deal with normal subgroups in the algebraic framework as we shall see later on. \\

\noindent
{\bf DEFINITION  2.13}: Two ideals $\mathfrak{a}, \mathfrak{b} \subset A$ are said to be {\it comaximal} if $\mathfrak{a} + \mathfrak{b} = A$.  \\

\noindent
{\bf LEMMA 2.14}: Let us consider ideals ${\mathfrak{a}}_1, ... , {\mathfrak{a}}_r \subset A$ that are comaximal two by two, that is ${\mathfrak{a}}_i + {\mathfrak{a}}_j = A$, $\forall i\neq j$. Then we have ${\mathfrak{a}}_1...{\mathfrak{a}}_r = {\mathfrak{a}}_1 \cap ... \cap {\mathfrak{a}}_r $.  \\

\noindent
{\it Proof}: For simplicity, we just consider the case $r=2$ with $\mathfrak{a}$ and $\mathfrak{b}$ two comaximal ideals in $A$, only proving that 
$\mathfrak{a} \mathfrak{b}=\mathfrak{a} \cap \mathfrak{b}$. The inclusion $\mathfrak{a}\mathfrak{b} \subseteq \mathfrak{a} \cap \mathfrak{b}$ being evident, we may find elements $a \in \mathfrak{a}, b \in \mathfrak{b}$  with $a + b = 1$. Hence, for any $x \in \mathfrak{a} \cap \mathfrak{b}$, we have $x= ax + bx$ with both $ax$ and $bx$ in $\mathfrak{a}\mathfrak{b}$, that is $x\in \mathfrak{a}\mathfrak{b}$. The general situation can be proved by induction ([28]).  \\
$\hspace*{12cm}  \Box $  \\

Let $A_1, ..., A_r$ be rings and consider all the sequences $\{ a_1, ..., a_r \mid a_i \in A_i \}$. We may provide to this set a structure of ring by setting:
\[  (a_1,..., a_r) + (b_1, ..., b_r) = (a_1 + b_1, ..., a_r + b_r), \,\,\, (a_1,...,a_r)(b_1,...,b_r) = (a_1b_1,...,a_r b_r)  \]
This ring is called the direct sum of $A_1,...,A_r$ and is denoted by $A_1 \oplus ... \oplus A_r$.  \\

\noindent
{\bf PROPOSITION  2.15}: ({\it Chinese remainder theorem}) With the same assumption as in the last Lemma, we have an isomorphism:  
\[    A / ({\mathfrak{a}}_1 \cap ... \cap {\mathfrak{a}}_r) \simeq (A/ {\mathfrak{a}}_1) \oplus ... \oplus (A/ {\mathfrak{a}}_r)   \]

An Artinian ring has only a finite number of prime ideals and each of them is maximal. As above, we shall denote them by ${\mathfrak{m}}_1, ... ,{\mathfrak{m}}_r $. \\

\noindent
{\bf THEOREM  2.16 }: Any Artinian ring with zero nilradical, that is with ${\mathfrak{m}}_1... {\mathfrak{m}}_r = 0$ is isomorphic to a direct sum of fields.  \\

\noindent
{\it Proof}: Any reduced residue ring of a polynomial ring in one indeterminate over a field $k$ is an Artinian ring. In fact, let $A=k[y]/ \mathfrak{a}$ where the principal 
ideal $\mathfrak{a}$ is generated by the polynomal $P \in k[y]$  that can be written $P=P_1...P_r$ with each $P_i \neq P_j $ irreducible in $k[y]$ over $k$ and 
relatively prime to $P_j, \forall i\neq j$. If ${\mathfrak{p}}_i$ is the prime and thus maximal principal ideal generated in $k[y]$ by $P_i$, then we have $ \mathfrak{a} = {\mathfrak{p}}_1 \cap ... \cap {\mathfrak{p}}_r $ and thus:
\[    A = k[y]/ \mathfrak{a} \simeq (k[y]/ {\mathfrak{p}}_1) \oplus ... \oplus (A/ {\mathfrak{p}}_r)  \]
Taking the residue ${\mathfrak{m}}_i$ of ${\mathfrak{p}}_i$  with respect to $\mathfrak{a}$, we obtain ${\mathfrak{m}}_1 ... {\mathfrak{m}}_r = 0$ and thus (make a picture):
\[   {\mathfrak{m}}_1 \cap ... \cap {\mathfrak{m}}_r = 0  \,\,\,  \Rightarrow \,\,\,A \simeq (A/ {\mathfrak{m}}_1) \oplus ... \oplus (A/ {\mathfrak{m}}_r)    \]
The commutative and exact diagram:  \\
\[ \begin{array}{rcccccccl}
  &  & 0  &  &  0  &  &  0  &  &   \\
  &  &  \downarrow & & \downarrow &  & \downarrow &  &    \\
0 & \rightarrow & \mathfrak{a} & \rightarrow & {\mathfrak{p}}_i & \rightarrow &  {\mathfrak{m}}_i & \rightarrow & 0 \\
   &  &  \parallel  & & \downarrow &  & \downarrow &  &   \\
0 & \rightarrow & \mathfrak{a} & \rightarrow & k[y] & \rightarrow & A  &  \rightarrow & 0   \\
  &  & \downarrow &  & \downarrow  & & \downarrow &  &  \\
  &   & 0  & \rightarrow & k[y] / {\mathfrak{p}}_i & \rightarrow  & A / {\mathfrak{m}}_i & \rightarrow & 0  \\
  &   &  &  & \downarrow & & \downarrow & &        \\
  &  &  &  &  0 &  & 0 & & 
\end{array}  \]
finally proves the isomorphism $A/ {\mathfrak{m}}_i \simeq k[y]/ {\mathfrak{p}}_i$, $\forall i= 1, ..., r$.  \\
$\hspace*{12cm}   \Box $ \\

\noindent
{\bf COROLLARY 2.17 }: $A \simeq  A_1 \oplus ... \oplus A_r  \,\,\, \Rightarrow \,\,\, Q(A) \simeq Q(A_1) \oplus ... \oplus Q(A_r)    $.  \\
 
\noindent
{\bf EXAMPLE 2.18 }: With $k= \mathbb{Q}, r=2$, let us consider the principal ideals $\mathfrak{a}, {\mathfrak{p}}_1, {\mathfrak{p}}_2$ in $k[y]$ generated respectively by $P=y^3 - 1, P_1= y - 1, P_2= y^2 + y + 1$. Denoting by $\eta$ the residue of $y$ in $A =k[y]/ \mathfrak{a}$, we have 
$\mathfrak{a} = {\mathfrak{p}}_1 \cap {\mathfrak{p}}_2$ and thus $A\simeq (k[y]/ {\mathfrak{p}}_1) \oplus  (k[y]/{\mathfrak{p}}_2)$.
We also obtain the Bezout identity:
\[   - \frac{1}{3} (y+2)(y-1) + \frac{1}{3} ( (y^2 + y + 1) =1   \]
Introducing the complex imaginary quantity $i$ with $i^2+1=0$ and the complex cube root of unity $j= (-1 + i \sqrt{3})/2$, we may set $A \simeq (\mathbb{Q} .1) \oplus (\mathbb{Q} 
 . 1 + \mathbb{Q} . j) \simeq \mathbb{Q} \oplus \mathbb{Q}(j) $ with $j^2 + j + 1=0 \Rightarrow j^3=1$. Accordingly, any element of $A$ may be written $ a= (\lambda, \mu + \nu j)$ with $\lambda, \mu, \nu \in \mathbb{Q}$.  \\

Recapitulating the results so far obtained, we have two procedures in order to compute the tensor product $K {\otimes}_k L$ of two fields containing a field $k$. \\

\noindent
1) Assuming for simplicity that the extension$K/k$ is finitely generated, we may exhibit a finite transcendence basis $(s)$ of $K/k$ in such a way that $K$ is algebraic over $k(s)$ while $k(s)$ is regular over $k$. This is the way that has been used in Theorem 2.10.   \\
2) We may also introduce the algebraic closure $K_0$ of $k$ in $K$, in such a way that $K_0$ is algebraic over $k$ and $K$ is regular over $K_0$. In this case, we have $K{\otimes}_k L = K {\otimes}_{K_0} (K_0 {\otimes}_k L)$. According to the previous Theorem, $K_0 {\otimes}_k L $ is isomorphic to a direct sum of fields, say 
$M_0 \oplus ... \oplus M_r$, each $M_i$ containing both $K_0$ and $L$. Thus, $K{\otimes}_k L \simeq  K {\otimes}_{K_0} (M_1 \oplus ... \oplus M_r) \simeq (K{\otimes}_{K_0} M_1) \oplus ... \oplus (K{\otimes}_{K_0} M_r)$. As $K/{K_0}$ is a regular extension, then each $K{\otimes}_{K_0} M_i$ is an integral domain and each $Q(K{\otimes}_{K_0} M_i)$ is a field.

Coming back to the classical Galois theory, let us consider an algebraic extension $L/K$ defined by $L= K[y]/ (P)$ where the polynomial $P$ is irreducible over $K$ and the principal ideal $(P)$ is thus prime in $K[y]$. We recall that $L/K$ is a {\it Galois extension} iff $L$ is a {\it splitting field} of $P$, that is $(P)$ decomposes over $L$ into the maximal principal ideals $(y - {\eta}_i)$ in $L[y]$ in which the ${\eta}_i \in L$ are the roots of $P$. As $L[y]/ (y - {\eta}_i) \simeq L, \forall i= 1,...,\mid L/K \mid$, we get: \\

\noindent
{\bf PROPOSITION  2.19 }: $L/K$ is a Galois extension if and only if $L {\otimes}_K L \simeq L \oplus ... \oplus L$ with $\mid L/K \mid$ terms or, equivalently, $L{\otimes}_K L \simeq L {\otimes}_{\mathbb{Q}} (\mathbb{Q} \oplus ... \oplus \mathbb{Q})$ with $\mid L/K \mid$ terms as $\mathbb{Q} \subseteq K \subset L$.   \\

\noindent
{\bf DEFINITION  2.20}: If $L/K$ s a finite algebraic extension, we shall denote by $split(L/K)$ {\it a} (care) smallest Galois extension of K containing $L$, which is defined up to an isomorphiosm over $L$.  \\

\noindent
{\bf EXAMPLE  2.21}: If $P=y^3 -2 \Rightarrow K=\mathbb{Q}\subset  L= \mathbb{Q} (\sqrt[3]{2}) \Rightarrow split(L/K)= L(j)$ with imaginary quantity $j= ( - 1 + i \sqrt{3})/2$ and $M=K(j)=K[z]/(z^2 + z +1)$ in the following diagram: \\
\[     \begin{array}{ccc}
L = \mathbb{Q}(\sqrt[3]{2}) & \stackrel{2}{\longrightarrow} & split(L/K)=(L,M)   \\
 3 \uparrow \hspace{2mm} & \nearrow 6 & \hspace{4mm} \uparrow  3   \\
K= \mathbb{Q}  & \underset 2{\longrightarrow }  & \hspace{3mm} \mathbb{Q}(j) = M
\end{array} \]
We notice that $L$ and $M$ are linearly disjoint over $K$ in $split(L/K)=(L,M)$. \\
On the contrary, setting $\eta=\sqrt[3]{2}$, $L=K(\eta)$ and $L'=K(j \eta)$, then $L$ and $L'$ are not linearly disjoint over $K$ in $split(L/K)=(L,L')$ because $(1, \eta, {\eta}^2)$ is a basis of $L$ over $K$ while $(1, j \eta, (j \eta)^2)$ is a basis of $L'$ over $K$ but we have the linear relation $ (j \eta)^2 \times 1 + (j \eta) \times \eta + 1 \times {\eta}^2 = 0 $ even though 
$L \cap L' =K$. As $\eta$ is a root of $y^3 - 2= 0$ while ${\eta}'= j \eta$ is a root of $(y')^3 - 2=0$ and ${\eta}' - \eta \neq 0$, we get $({\eta}')^2 + \eta {\eta}' + {\eta}^2 =0$ because 
$j^2 + j + 1 = 0$.\\

More generally, if $\eta$ is a generic zero of a prime ideal $\mathfrak{p}$ such that $L =K(\eta)=Q(K[y]/ \mathfrak{p})$, then the perfect ideal $M \mathfrak{p} \subset M[y]$ is a finite intersection ${\mathfrak{p}}_1 \cap ... \cap {\mathfrak{p}}_r$ of prime ideals and we may define $Q(L {\otimes}_K M)= Q(M[y]/ M \mathfrak{p})$ as a direct sum of fields. The situation of 
$L {\otimes}_K L$ is slightly different as we saw it is a reduced Artinian ring having only a finite number of of prime ideals ${\mathfrak{m}}_1, ..., {\mathfrak{m}}_r$ that are also maximal with zero intersection and we have the direct sum of fields: 
  \[  Q(L {\otimes}_K L) \simeq Q(L {\otimes}_K L/ {\mathfrak{m}}_1) \oplus ... \oplus Q( L {\otimes}_K  L/ {\mathfrak{m}}_r)  \]
  
As we do not need to take the full rings of quotients, using the commutative diagram:
\[ \begin{array}{ccccccccc}
& & & & 0 & & 0 & &   \\
& & & & \downarrow & & \downarrow & & \\
& & 0 & \rightarrow    & L & = & L & \rightarrow & 0  \\
& & \downarrow & & \downarrow & & \downarrow & & \\
0 & \longrightarrow & {\mathfrak{m}}_i & \longrightarrow & L{\otimes}_K L & \longrightarrow & L {\otimes}_K L / {\mathfrak{m}}_i & \longrightarrow  & 0
\end{array}   \]
we shall call ${\eta}_i$ the image of $\eta$ under the composite monomorphism: 
\[  L \rightarrow 1 \otimes L \rightarrow L {\otimes}_K L\rightarrow   L {\otimes}_K L/ {\mathfrak{m}}_i  \]
set $L_i=K({\eta}_i)$ and identity $ \eta$ with its image under the composite monomorphism:
\[  L \rightarrow L \otimes 1 \rightarrow L {\otimes}_K L\rightarrow   L {\otimes}_K L/ {\mathfrak{m}}_i  \]
In actual ractice, if $P$ is the minimal unitary polynomial in $K[y]$ of a primitive element of $L/K$, we may denote by $\bar{P}$ the image of $P$ under the isomorphism $K[y] \rightarrow K[\bar{y}]:y \rightarrow \bar{y}$ in order to obtain $Q(L{\otimes}_K L)\simeq Q(K[y,\bar{y}]/(P,\bar{P}))$. As $\bar{P} - P$ is divisible by $(\bar{y} - y)$, we may label the ${\mathfrak{m}}_i$ in such a way that ${\eta}_1= \eta \Rightarrow L_1= K({\eta}_1)=K(\eta)=L$.\\
We obtain $Q(L{\otimes}_K L / {\mathfrak{m}}_i) = K(\eta, {\eta}_i) = L({\eta}_i) = L_i (\eta) = (L,L_i) = N_i $ and the useful formula: 
\[  {\Sigma}_i \mid (L,L_i) /L \mid = \mid (L,L_1) / L \mid + ... + \mid (L,L_r)/L \mid= \mid L/K \mid   \]
We finally define a {\it finite family} ${\sigma}_1, ..., {\sigma}_r \in iso(L/K)$ with ${\sigma}_1=id_L$ by the formula ${\sigma}_i(\eta]={\eta}_i$ in such a way that $\sigma \in iso (L/K) \Rightarrow K (\sigma (\eta) \simeq L_i$ for some $i=1,...,r$. This definition only depends on the finite extension $L/K$. \\

\noindent
{\bf DEFINITION  2.22}: We shall say that an isomorphism $\sigma$ is {\it conjugate} of ${\sigma}_i$ if $\sigma (\eta)$ and ${\eta}_i$ are roots of the same minimum polynomial over $L$. When dealing with finite extensions, any $\sigma \in iso (L/K)$ is the conjugate of one and only one ${\sigma}_i$ as already defined. For this reason, we shall say that each isomorphisms 
${\sigma}_i$ is an {\it isolated isomorphism}. \\

\noindent
{\bf LEMMA  2.23 }: One has $inv({\sigma}_1, ..., {\sigma}_r) = K$. \\ 

\noindent
{\it Proof}: Let $\eta$ be a primitive element of $L/K$ with a minimum unitary polynomial $P$ of degree $m$. Accordingly, every element of $L$ can be written as a polynomial in $\eta$ of degree $\leq (m-1)$. Let thus $\zeta = a_1 {\eta}^{m-1} + ... + a_m$ with $a_1, ..., a_m \in K$ be such an element satisfying ${\sigma}_i(\zeta) = \zeta, \forall i= 1,...,r$. As we are dealing with principal ideals, we have the prime decomposition $L (P) = (P_1) \cap ... \cap (P_r)$. Now, the polynomial $R(y) = a_1 y^{m-1} + ... + a_{m-1} y - (a_1 {\eta}^{m-1} + ... + a_{m-1} \eta) \in L[y]$ is such that $R({\eta}_1)= ... = R({\eta}_r)=0$ by assumption. According to Hilbert theorem we must have $R \in L(P)$ and this is impossible unless 
$a_1 = ... = a_{m-1} = 0$, a result leading to $\zeta = a_m \in K$.  \\
\hspace* {12cm}   $\Box $  \\

\noindent
{\bf DEFINITION  2.24 }: We say that $L/K$ is an {\it automorphic extension} if a model variety $\Sigma$ of $L/K$ is a {\it principal homogeneous space} (PHS) for a finite algebraic group 
$\Gamma = \Gamma (L/K)$ defined over $k \subset K$ in such a way that each component of $k[\Gamma]$ is linearly disjoint from $L$ over $k$. We have the {\it fundamental isomorphism} $Q(L {\otimes }_K  L)  \simeq  Q (L {\otimes}_k  k[\Gamma])$ but the full rings of quotient may not be needed. For simplicity, we shall suppose that $L/k$ is a regular 
extension in order to be more coherent with the point of view adopted by Vessiot and we shall say that $L/K$ is {it regular} over $k$. \\

\noindent
{\bf EXAMPLE  2.25 }: If we consider the finite extension $L/K$ with $K=\mathbb{Q} \subset L=\mathbb{Q}(\eta)$ with $\eta=\sqrt[8]{2})$ the (real) generic zero of the underlying irreducible equation is $y^8 - 2 = 0$, we have thus  $\sqrt{2} = {\eta}^4 \subset L$ and we obtain successively (See any textbook for {\it cyclotomic fields}):  \\
\[  \begin{array}{rcl}
{\bar{y}}^8 - y^8 & =  & (\bar{y} - y)(\bar{y} + y)({\bar{y}}^2 + y^2)({\bar{y}}^4 + y^4)   \\
   &   =  & (\bar{y} - y)(\bar{y} + y)( \bar{y} - i y)(\bar{y} + i y) (\bar{y} - \frac{1+ i}{\sqrt{2}} y) (\bar{y} + \frac{1+i}{\sqrt{2}} y) (\bar{y} - \frac{1-i}{\sqrt{2}}  y)
(\bar{y} + \frac{1-i}{\sqrt{2}} y)    
\end{array}  \] 
Setting $\alpha=(1+i)/{\sqrt{2}} \Rightarrow {\alpha}^2 = i \Rightarrow  \sqrt{2} = \alpha + \frac{1}{\alpha}$, we notice that $\mathbb{Q}(\alpha) = Q(\mathbb{Q}[a]/ (a^4 + 1))$ and $L$ are {\it not} linearly disjoint in 
$split(L/K)=L(\alpha)$ over $ \mathbb{Q}$ because we have $1 \times ({\alpha}^2 + 1) - \sqrt{2} \times \alpha  = 0$ for bases $(1, \eta, {\eta}^2, ..., {\eta}^4=\sqrt{2}, ..., {\eta}^7)$ for $L/K$ and $(1, \alpha, {\alpha}^2, {\alpha}^3)$ for $\mathbb{Q}(\alpha)/ \mathbb{Q}$. \\
On the contrary, with $k=\mathbb{Q}\subset K=k(y^8) \subset L=k(y)$, then $\mathbb{Q}(\alpha)$ and $L$ are linearly disjoint in $L(\alpha)$ over $k$. In the present situation, the group of invariance is the cyclic group $\Gamma$ with generator $\bar{y} = a y$ and $k[\Gamma]= k[a]/(a^8- 1)=k \oplus  k \oplus  k[a]/(a^2 + 1) \oplus  k[a]/ (a^4 + 1)=M_1 \oplus M_2 \oplus M_3 \oplus M_4$ with $M_1= M_2 = k, M_3 = k(i), M_4 = k (\alpha)$.  \\
We finally consider the irreducible equation $y^{12} - 3=0$ over $\mathbb{Q}$ by using the fact that:  \\
\[   {\bar{y}}^{12} - y^{12} = (\bar{y} - y)(\bar{y} + y)({\bar{y}}^2 + y {\bar{y}} + y^2)({\bar{y}}^2 - y {\bar{y}} + y^2)({\bar{y}}^2 + y^2)({\bar{y}}^4 - y^2 {\bar{y}}^2 + y^4)  \]
It follows that $\mathbb{Q}(\eta)$ with $\eta = \sqrt[12]{3} $ is linearly disjoint over both $M_3 = \mathbb{Q}(j)$ and $M_5 = \mathbb{Q}(i)$ or even $M_6 = \mathbb{Q}(ij)$ but {\it not} over $\mathbb{Q}(i, j)$ as this later field contains ${\eta}^6 = \sqrt{3} = (1 + 2 j)/ i $.  \\

Developing the tensor products that appear in the fundamental isomorphism, we obtain by linear disjointness:\\

\noindent
{\bf PROPOSITION  2.26 }: One has $ Q( L {\otimes}_K L ) / {\mathfrak{m}}_i = (L,L_i) = (L, M_i) = Q( L{\otimes}_k M_i) $, $\forall i = 1, ... , r$:
\[    \begin{array}{rcl}
 Q( L {\otimes}_K L )& \simeq  & (L,L_1) \oplus ... \oplus (L,L_r)   \\
                                &  \simeq  &  Q(L {\otimes}_k M_1) \oplus ... \oplus  Q(L {\otimes}_k M_r) \\
                                &  \simeq  &  Q(L {\otimes }_k ( M_1 \oplus ... \oplus M_r))   \\
                                &   \simeq  &  Q(L { \otimes}_k k[\Gamma])
\end{array}   \]
and the commutative diagram:   \\

\[     \begin{array}{cccccc}
L &  \longrightarrow  & N_i  &  \longrightarrow  &  split(L/K) \\
\uparrow  &    &  \uparrow    &    &           \uparrow     \\
K  & \longrightarrow      & (K,M_i) &    \longrightarrow      & (K,M)  \\
\uparrow  &    &  \uparrow    &    &  \uparrow                \\
k  & \longrightarrow      & M_i  &    \longrightarrow    & M
\end{array}  \]

\noindent
in which $M = (M_1, ... , M_r)$ and $split(L/K) = (L, M) = (L_1, L_2, ... , L_r)$ because $L_1=L$.  \\

From now on, in order to simplify the proofs while showing the usefulness of these methods as we shall only use fields, we shall suppress the rings of quotients $Q(  \,\, )$. \\

\noindent
{\bf PROPOSITION  2.27}: $k[\Gamma] \subset L {\otimes}_K L$ is the ring of polynomial functions on an algebraic finite group $\Gamma = \Gamma (L/K)$ defined over $k$ and has thus 
an induced structure of {\it Hopf algebra} because $k \subset k[\Gamma]$.  \\

\noindent
{\it Proof}: Let us introduce the two monomorphisms:
\[  \left\{  \begin{array}{lll}
 \tilde{\alpha} & :  &  L \longrightarrow  L {\otimes}_K L : a \rightarrow a \otimes 1  \\
 \tilde{\beta}  & :   &  L \longrightarrow  L {\otimes}_K L : b \rightarrow 1 \otimes  b
\end{array} \right.  \]
respectively called {\it source inclusion} and {\it target inclusion} with the exact double arrow sequence:  
\[ 0 \longrightarrow  K  \longrightarrow L  \underset{\tilde{\beta}}{\stackrel{\tilde{\alpha}}{\rightrightarrows}} L {\otimes}_K L  \]

Now, we have the isomorphisms:
\[ L {\otimes}_K L {\otimes}_K L \simeq L {\otimes}_K L {\otimes}_k k[\Gamma] \simeq L {\otimes}_k k[\Gamma] {\otimes}_k k[\Gamma] \]
and the top row of the commutative diagram:
\[ \begin{array}{ccccc}
0 & \longrightarrow & L {\otimes}_k k[\Gamma]  & \longrightarrow & L {\otimes}_k k[\Gamma] {\otimes}_k k[\Gamma]   \\
  & & \downarrow   &  & \downarrow \\
0 & \longrightarrow & L {\otimes}_K L & \longrightarrow  & L {\otimes}_K L {\otimes}_K L  \\
 &  &  &  &  \\
   &  &   a \otimes b  & \longrightarrow  & a \otimes 1 \otimes b 
\end{array}  \] 
induces the {\it diagonal comorphism} $\tilde{\epsilon}: k[\Gamma] \rightarrow k[\Gamma] {\otimes}_k k[\Gamma]$ by linear disjointness over $k$. In actual practice, we have for example 
$\bar{y} = a y, \bar{\bar{y}}=b \bar{y}\Rightarrow \bar{\bar{y}} = b a y $ for the multiplicative group, that is $a=\frac{\bar{y}}{y}= \frac{1}{y} \otimes y \otimes 1, b=\frac{\bar{\bar{y}}}{\bar{y}} = 1 \otimes \frac{1}{y} \otimes y \Rightarrow ba =  \frac{\bar{\bar{y}}}{\bar{y}}= \frac{1}{y} \otimes 1 \otimes y $ (Compare to [3] Th 1, p 97). \\
Let us study with more details this monomorphism. Substituting, we get:

 \[ k[\Gamma] = {\oplus}_i M_i \rightarrow  k[\Gamma] {\otimes}_k k[\Gamma] = ({\oplus}_i M_i ) ({\oplus}_j M_j) = {\oplus}_{i,j} (M_i {\otimes}_k M_j) \]

Then, the top row of the following commutative diagram:
\[   \begin{array}{ccccl}
L {\otimes}_k k[\Gamma] & \rightarrow & L & \rightarrow  &  0  \\
\downarrow &   &\parallel &  &   \\
L {\otimes}_K L &  \rightarrow & L & \rightarrow & 0 \\
  &   &   &  &   \\
a \otimes b & \rightarrow  &  ab &   &
\end{array}  \]
induces the {\it augmentation comorphism} $\tilde{id}: k[\Gamma] \rightarrow k \rightarrow k[\Gamma]$ by linear disjointness over $k$. In actual practice, we have 
$ \bar{y} = a y \Rightarrow  a = \frac{\bar{y}}{y}= \frac{1}{y} \otimes y \Rightarrow  \frac{1}{y} \otimes  y \rightarrow  1$ for the identity $\bar{y} = 1 y = y$.  \\

Finally, the top row of the commutative diagram:
\[   \begin{array}{ccccl}
L {\otimes}_k k[\Gamma] & \rightarrow & L {\otimes}_k k[\Gamma] & \rightarrow  &  0  \\
\downarrow &   &\downarrow &  &   \\
L {\otimes}_K L &  \rightarrow & L {\otimes}_K L & \rightarrow & 0 \\
  &   &   &  &   \\
a \otimes b & \rightarrow  & b \otimes a&   &
\end{array}  \]
induces the {\it antipode comorphism} $\tilde{\iota} : k[\Gamma] \rightarrow k[\Gamma] $ by linear disjointness over $k$, which  is sending each $M_i$ to itself as can be easily seen on each preceding example. In actual practice, we have $\bar{y}=a y \Rightarrow  a = \frac{\bar{y}}{y} = \frac{1}{y} \otimes y \Rightarrow {a}^{-1}= \frac{1}{a} = \frac{y}{\bar{y}}= y \otimes \frac{1}{y}$. \\
\hspace*{12cm}  $\Box$ \\

\noindent
{\bf THEOREM 2.28}: If $L/K$ is an automorphic extension regular over $k$ for a group $\Gamma$ and $K \subset K' \subset L$ is an intermediate field, then $L/K'$ is an automorphic extension regular over $k$ for a subgroup ${\Gamma}' \subset \Gamma$.  \\

\noindent
{\it Proof}: As $iso(L/K') \subset iso(L/K)$, we may use a convenient labelling such that the isomorphisms ${\sigma}_1, ... ,{\sigma}_s \in iso(L/K') $ constructed as before are among the isomorphisms ${\sigma}_1, ..., {\sigma}_r \in iso(L/K) $ with $s < r$ in such a way that ${\sigma}_1 = id_L$. In actual practice, if $\mathfrak{r}$ is the ideal of $L {\otimes}_K L$ generated by all the elements of the form $a\otimes 1 - 1 \otimes a $ with $a \in K'$, then $L {\otimes}_K L / \mathfrak{r} \simeq L {\otimes }_{K'} L$. If now ${\mathfrak{m}}'_i$ is a prime ideal of 
$L {\otimes}_{K'} L$, the {\it inverse image} of ${\mathfrak{m}}'_i$ under the canonical epimorphism $L {\otimes}_K  L \rightarrow L {\otimes}_{K'}  L$ is a prime ideal of $L {\otimes }_K  L$ and must therefore be equal to some ${\mathfrak{m}}_i$. Using a chase in the following commutative and exact diagram:  \\
\[  \begin{array}{ccccccccc}
  &  & 0 & & 0 & & & & \\
   &    & \downarrow  & & \downarrow   & & & &  \\
   0 & \longrightarrow & \mathfrak{r} & =  &  \mathfrak{r}   &  \longrightarrow & 0  & & \\
   & & \downarrow    & & \downarrow  &  & \downarrow &  & \\
   0 & \longrightarrow & {\mathfrak{m}}_i & \longrightarrow & L {\otimes}_K L & \longrightarrow & (L{\otimes}_K L)/{\mathfrak{m}}_i & \longrightarrow  & 0  \\
   & & \downarrow    & & \downarrow  &  & \downarrow &  & \\
   0 & \longrightarrow & {\mathfrak{m}}'_i & \longrightarrow & L {\otimes}_{K'} L & \longrightarrow & (L{\otimes}_{K'} L)/{\mathfrak{m}}'_i & \longrightarrow  & 0  \\
   & & \downarrow & & \downarrow  & & \downarrow   & & \\
   & & 0 & & 0 & & 0 & &   
   \end{array}    \]
we obtain an isomorphism $L {\otimes}_{K'} L / {\mathfrak{m}}'_i \simeq L {\otimes}_K L / {\mathfrak{m}}_i = N_i$. \\
Now, according to Proposition 2.26, we have $N_i = Q(L{\otimes }_k M_i) = (L,M_i) \Rightarrow  M_i = N_i \cap M$. In addition, we have 
$s < r \rightarrow split (L/K') \subseteq split (L/K)$ in the following commutative diagram:  \\

\[     \begin{array}{ccccccc}
L &  \longrightarrow  & N_i  &  \longrightarrow  &    split(L/K') & \longrightarrow & split(L/K) \\
\uparrow  &    &  \uparrow    &    &     \uparrow     &       &     \uparrow     \\
K'   & \longrightarrow      & (K',M_i ) &    \longrightarrow   & (K', M' ) &  \longrightarrow  &   (K',M)  \\
\uparrow  &    &  \uparrow    &    &  \uparrow       &       &     \uparrow     \\
K  & \longrightarrow      & (K,M_i) &    \longrightarrow   &  \uparrow &  \longrightarrow  & (K,M)  \\
\uparrow  &    &  \uparrow    &    &  \uparrow       &       &     \uparrow     \\
k  & \longrightarrow      & M_i  &    \longrightarrow   &  M'&  \longrightarrow  & M
\end{array}  \]
with $M' = (M_1,..., M_s) \subset (M_1, ..., M_r) = M$ and still $L$ linearly disjoint over $M'$ in $split(L/K')$.\\
It remains to prove that ${ \Gamma}' $ can be constructed like $\Gamma $ through its Hodge algebra as in ([3] Lemma 1, p 93 and Lemma 5, p 101). For this, introducing the ideal 
$ \mathfrak{s}= \mathfrak{r} \cap k[\Gamma] \subset L {\otimes}_k L$ in the following commutative and exact diagram:

\[  \begin{array}{rcccc}
  &   & 0  &   &  0   \\
  &    &  \downarrow &  &  \downarrow  \\
0  &  \longrightarrow  &  \mathfrak{s} &   \longrightarrow  &  \mathfrak{r}  \\
  &   &  \downarrow  &  &  \downarrow  \\
0 & \longrightarrow & k[\Gamma ] & \longrightarrow & L {\otimes}_K L  \\
   &   &     \downarrow &   &  \downarrow   \\
0 & \longrightarrow & k[ {\Gamma}' ] & \longrightarrow & L {\otimes}_{K'} L  \\
   &                         &  \downarrow &       &  \downarrow        \\
   &             &  0  &    &   0     
   \end{array}  \]
we obtain the {\it specialization} epimorphism $k[\Gamma] \rightarrow k[{\Gamma}'] \rightarrow 0$. If ${\mathfrak{n}}_i = {\mathfrak{m}}_i \cap k[\Gamma]$ and 
${\mathfrak{n}}'_i = {\mathfrak{m}}'_i \cap k[{\Gamma}']$, then ${\mathfrak{n}}_i $ is the inverse image of ${\mathfrak{n}}'_i$  by this epimorphism. A chase in the following commutative and exact diagram:  \\
\[  \begin{array}{ccccccccc}
  &  & 0 & & 0 & & & & \\
   &    & \downarrow  & & \downarrow   & & & &  \\
   0 & \longrightarrow & \mathfrak{s} & =  &  \mathfrak{s}   &  \longrightarrow & 0  & & \\
   & & \downarrow    & & \downarrow  &  & \downarrow &  & \\
   0 & \longrightarrow & {\mathfrak{n}}_i & \longrightarrow &k[\Gamma] & \longrightarrow & k[\Gamma]/{\mathfrak{n}}_i & \longrightarrow  & 0  \\
   & & \downarrow    & & \downarrow  &  & \downarrow &  & \\
   0 & \longrightarrow & {\mathfrak{n}}'_i & \longrightarrow & k[{\Gamma}'] & \longrightarrow & k[{\Gamma}']/{\mathfrak{n}}'_i & \longrightarrow  & 0  \\
   & & \downarrow & & \downarrow  & & \downarrow   & & \\
   & & 0 & & 0 & & 0 & &    
   \end{array}    \]
finally proves the isomorphism $k[{\Gamma}'] / {\mathfrak{n}}'_i \simeq k[\Gamma] / {\mathfrak{n}}_i = M_i$ and we have thus $k[{\Gamma}'] \simeq M_1 \oplus ... \oplus M_s$. We obtain therefore: 
\[   \begin{array}{ccl}
L {\otimes}_k k[{\Gamma}'] &  \simeq  &  L {\otimes}_k (M_1\oplus ... \oplus M_s)  \\
                                           &  \simeq  &  (L {\otimes}_k M_1) \oplus ... \oplus (L {\otimes}_k M_s)  \\
                                          &   \simeq  &  N_1 \oplus ... \oplus N_s  \\
                                          &  \simeq &  L {\otimes}_{K'} L
  \end{array}  \]
and the following commutative and exact diagram:
\[  \begin{array}{ccccccc}
0  &  \longrightarrow & L {\otimes}_k k[\Gamma] & \longrightarrow &  L {\otimes}_K L  & \longrightarrow & 0  \\
    &                           & \downarrow &       &  \downarrow   &   & \\
0  &  \longrightarrow & L {\otimes}_k k[{\Gamma}'] & \longrightarrow &  L {\otimes}_{K'} L  & \longrightarrow & 0  \\
    &   &     \downarrow  &    &  \downarrow   &   &    \\
    &   &    0   &   &  0  &  &
\end{array}   \]
proving that $L/K'$ is an automorphic extension for ${\Gamma}' = \Gamma (L/K') \subset \Gamma(L/K)=\Gamma$.  \\
\hspace*{12cm}  $\Box $ \\

Let us now consider an automorphic extension $L/K$ regular over $k$ for an algebraic finite group $\Gamma = \Gamma (L/K)$ defined over $k$. We have proved that, if $ K \subset K' \subset L$ for an intermediate field $K'$, then  $L/K'$ is an automorphic extension for an algebraic subgroup ${\Gamma}' = \Gamma(L/K') \subset \Gamma$ defined over $k$. The problem left is thus to study the normality condition ${\Gamma}' \lhd \Gamma $ by finding a criterion involving only the three field extensions $K, K', L$ of $k$ but such a result is not intuitive at all.  \\

\noindent
{\bf DEFINITION  2.29}: The composite {\it translation comorphism} $\tilde{\tau}: L \stackrel{\tilde{\beta}}{ \longrightarrow} L {\otimes}_K L \longrightarrow L {\otimes}_k k[\Gamma]$ is the comorphism of a rational action of  $\Gamma$ on a model variety $\Sigma$ of $L/K$.  \\

\noindent
{\bf LEMMA  2.30}: One has ${\Gamma}' \lhd \Gamma$ {\it if and only if} the action of $\Gamma$ on $L/K$ over $k$ induces an action of $\Gamma$ on $K'/K$ over $k$ acording to the following commutative and exact diagram for the translation comorphism $\tilde{\tau}$: \\
\[ \begin{array}{ccccc}
0  &  &  0 &  & 0   \\
\downarrow &  &  \downarrow    & \hspace{2cm} & \downarrow   \\
  K' & \stackrel{\tilde{\tau}}{\longrightarrow} & K' {\otimes}_k k[\Gamma] &\hspace{2cm} & K(\Omega) = K'   \\
 \downarrow & & \downarrow &  & \downarrow  \\
 L & \stackrel{\tilde{\tau}}{\longrightarrow} & L {\otimes}_k k[\Gamma]  & \hspace{2cm}  & K(\Sigma) = L 
 \end{array}      \]

\noindent
{\it Proof}: If ${\Gamma}' \lhd \Gamma$ and $\eta, {\eta}'$ are two generic points of $\Sigma$ such that ${\eta}' = h \eta$ for a certain $h \in {\Gamma}'$ and $g\in \Gamma$, then $\exists h' \in {\Gamma}'$ such that $ g {\eta}' =g h \eta = h' g \eta$. Conversely, if the action of $\Gamma$ passes to the quotient  on $\Omega = \Sigma / {\Gamma}' $, then $g {\eta}' = g h \eta = h' g \eta $ for a certain $h' \in {\Gamma}'$ and thus $ g h = h' g $ for any $g \in \Gamma$, because the action of $\Gamma $ is free, that is to say ${\Gamma}' \lhd \Gamma$. We have the picture: 
\[   \begin{array}{ccccc}
 & & \eta & \stackrel{g}{\longrightarrow} & g \eta \\
 \Sigma &\hspace{2cm} & h \downarrow \,\,\, &  & \,\,\, \downarrow  h'   \\
  &    & {\eta}' &\stackrel{g}{\longrightarrow} & g {\eta}' \\
   &  &   \vdots       &         &  \vdots   \\
  \Omega &  \hspace{2cm} &  \bullet & \stackrel{g}{ \longrightarrow}  & .\bullet
  \end{array}. \]
It is finally sufficient to notice that $\Omega$ is a model variety for $K'/K$.  \\

In a more practical way, the chain of inclusions $K \subset K' \subset L $ provides a chain of inclusions $ K' {\otimes}_K K' \subset L {\otimes}_K K' \subset L {\otimes}_K L$. Accordingly, any prime and thus maximal ideal ${\mathfrak{m}}'_i \subset K' {\otimes}_K K'$ can be extended to a perfect ideal of $L {\otimes}_K K'$ because $K' \subset L$, which is an intersection of prime ideals in $L {\otimes}_K K'$ and each such prime ideal can be similarly extended to a perfect ideal in $L {\otimes}_K L$. Hence, each ${\mathfrak{m}}''_i$ can be extended to a certain prime and thus maximal ideam ${\mathfrak{m}}_i \subset L { \otimes}_K L$. It follows that each isolated isomorphism $ {\sigma}'_i \in iso (K'/K)$ can be extended to an isolated isomorphism ${\sigma}_i \in iso (L/K)$.  \\
\hspace*{12cm}      $ \Box $. \\
 
\noindent
{\bf DEFINITION 2.31}: If $L/K$ is an automorphic extension over a field $k$, an intermediate automorphic extension $K'/K$ is said to be 
{\it admissible} if one has the following commutative and exact diagram of reciprocal image for $G= \Gamma (K'/K)$: 
\[   \begin{array}{ccccc}
  &     & 0 & & 0 \\
  & & \downarrow & & \downarrow \\
  0 & \longrightarrow & k[G] & \longrightarrow & K' {\otimes}_K K'  \\
   & & \downarrow &  & \downarrow \\
   0 & \longrightarrow & k[\Gamma] & \rightarrow & L {\otimes}_K L 
   \end{array}. \]

If $L/K$ is an automorphic extension over a field $k$  for an algebraic group $\Gamma = \Gamma (L/K)$, we have proved that there exists a bijective dual correspondence between the algebraic subgroups of $\Gamma$ defined over $k$ and the intermediate fields between $K$ and $L$. The following example will prove that $K'/K$ may be an automorphic extension even if ${\Gamma}'$ is {\it not} normal in $\Gamma$, when $K'/K$ is {\it not} admissible, contrary o the classical Galois theory.   \\

\noindent
{\bf EXAMPLE 2.32}: With $k=K= \mathbb{Q}, K' = \sqrt[3]{2}= \mathbb{Q}(\eta), L= \mathbb{Q} (\sqrt[3]{2}), j)$ with $j^2 + j + 1=0 \Rightarrow j^3=1$, we already know that $\mid L/K \mid = \mid L/K' \mid \times \mid K'/K \mid = 2 \times 3 = 6$ and that $L/K$ is a Galois extension for $S_3$, thus an automorphic extension as in Example 1 with ${\omega}^1= {\omega}^2 = 0, {\omega}^3 = 2$. We may define $\sigma (\eta)=j \eta, {\sigma}^2 (\eta)=j^2 \eta, \sigma (j)=j$ and $\tau(\eta)=\eta, \tau(j)= j^2$, that is $\Gamma = \{ e, \sigma, \tau \}$ has $3$ generators but $\mid \Gamma \mid = 6$ indeed with $\Gamma= \{e, \sigma, {\sigma}^2, \tau, \sigma \tau, \tau \}$ with $\sigma \tau \neq \tau \sigma$ as $\tau \sigma (\eta) = j^2 \eta$ while $\sigma \tau (\eta) = j^2 \eta = {\sigma}^2 \tau (\eta)$ and we finally notice that ${\sigma}^3= e, {\tau}^2=e$. However, we also know that $K'/K$ is {\it not} a Galois extension because ${\Gamma}' = \{ e, \tau\}$ is not normal in $\Gamma$ as we have ${\sigma}^{-1}\tau \sigma  = {\sigma}^2 \tau \sigma= {\sigma}^4 \tau = \sigma \tau \notin \Gamma$. We have finally $k[\Gamma] \simeq \mathbb{Q} \oplus ... \oplus \mathbb{Q}$ with $6$ terms while $k[G]\simeq \mathbb{Q} \oplus \mathbb{Q}(j)$ that is we have $K' {\otimes}_K K' \subset L {\otimes}_K L$ indeed but $k[G] \notin k[\Gamma]$. {\it On the contrary}, if we choose $K'= K(j)$ which s a Galois extension of $K$ with $\mid K'/K \mid = 2$ both with $L/K'$ which is a alois extension for $A_3 \lhd S_3$ with $\mid L/K' \mid = 3$ and ${\Gamma}' = \{ e, \sigma, {\sigma}^2 \}$ We have {\it now} $K' {\times}_K K' $ defined by $z^2 + z + 1=0, {\bar{z}}^2 + \bar{z} + 1 =0 \Rightarrow (\bar{z} - z)(\bar{z} + z + 1)$, a result leading to $k[G] \simeq \mathbb{Q} \oplus \mathbb{Q} \subset k[\Gamma]$ because $\bar{z}=z=j$ or $\bar{z}= -(z+1)= -(j+1) = j^2$ or, equivalently, $\tau (j)=j, {\tau}^2 (j)=j^2 = -(j+1)$ for $G \simeq S_3 / A_3 = \{ e, \tau \}$. Such a example explains why only admissible extensions $K'/K$ must be considered when normality is involved.  In this case, we have:   \\

\noindent
{\bf THEOREM 2.33}: An intermediate field $K'$ is an automorphic extension of $K$ for a group $G$ defined over $k$ if and only if ${\Gamma}' = \Gamma (L/K')$is a normal subgroup of $\Gamma$.  \\

\noindent
{\it Proof}: Composing with the morphism $\tilde{\beta}$, we obtain the commutative composite diagram: \\
\[   \begin{array}{ccccc}
K' & \stackrel{\tilde{\beta}}{\longrightarrow} & K' {\otimes}_K K' & \longrightarrow & K' {\otimes}_k k[G]. \\
\downarrow &  & \downarrow &    & \downarrow  \\
L & \stackrel{\tilde{\beta}}{\longrightarrow} & L {\otimes}_K L & \longrightarrow  & L {\otimes}_k k[\Gamma]
\end{array}      \]
in which the vertical arrows are monomorphisms and the inclusion $k[G] \subset k[\Gamma]$ is induced by the inclusion $K' {\otimes}_K K' \subset L {\otimes}_K L$ while the specialization $k[\Gamma] \rightarrow k[{\Gamma}'] \rightarrow 0$ is induced by the specialization $L{\otimes}_K L \rightarrow L {\otimes}_{K'} L \rightarrow 0$. According to the preceding lemma, we obtain therefore ${\Gamma}' \lhd \Gamma$ and $\Gamma / {\Gamma}' \simeq G$.  \\
{\it Conversely}, as $k[\Gamma] \subset L {\otimes}_K L$, introducing the Hopf algebra $k[G] = (K' {\otimes}_K K' )\cap k[\Gamma] \subset L {\otimes}_K L$ like in the last definition,  we have the commutative composite diagram:
\[    \begin{array}{ccccc}
K' {\otimes}_k k[G] & \longrightarrow & K' {\otimes}_K K' & \longrightarrow & K' {\otimes}_k k[\Gamma]     \\
\downarrow & & \downarrow & lemma & \downarrow \\
L {\otimes}_k k[\Gamma] & \longrightarrow & L {\otimes}_K L & \longrightarrow & L {\otimes}_k k[\Gamma]
\end{array}    \]
in which the three vertical arrows are monomorphisms.    \\
By composition, the morphism $K' {\otimes}_k k[G] \rightarrow K' {\otimes}_K K' $ is thus a monomorphism, because $k[G] \subset k[\Gamma]$. However, as $K' \subset L$ and $k[G] \subset   k[\Gamma] $, the left  diagram is just an inverse image and the later morphism  is also an epimorphism because $L/K$ is an automorphic  extension. Needless to say that, when $L/K, L/K', K'/K$ are Galois extensions, this well known result of normality becomes evident because $L {\otimes}_K L$ is a direct sum of $\mid L/K \mid$ copies of  $L$, $L {\otimes}_{K'} L$ is a direct sum of $\mid L/K' \mid $ copies of $L$ while $K' {\otimes}_K K'$ is a direct sum of $\mid K'/K \mid$ copies of $K'$ with $\mid L/K \mid = \mid L/K' \mid \times  \mid K'/K \mid $, a result leading to $\mid \Gamma \mid = \mid {\Gamma}' \mid \times \mid G \mid $.   \\                                                   
\hspace*{12cm}. $ \Box $ \\

\noindent
{\bf REMARK 2.34}: Though we have an epimorphism $k[\Gamma] \rightarrow  k[ {\Gamma}'] \rightarrow 0$ and a monomorphism $0 \rightarrow k[G] \rightarrow k[\Gamma] $, the later does not define the kernel of the previous one, as can be seen in the last example by using the commutative diagram:
\[  \begin{array}{ccc}
K' {\otimes}_K K' & \longrightarrow & K'  \\
\downarrow &  & \downarrow   \\
L {\otimes}_K L & \longrightarrow &  L {\otimes}_{K'} L 
\end{array} \]
in which the left vertical arrow is a monomorphism, the upper arrow is the epimorphism $a \otimes b \rightarrow ab$ while the right vertical arrow is the monomorphism $ a \rightarrow a \otimes 1 = 1 \otimes a$. In any case we have the recapitulating diagram:  \\
\[  \begin{array}{ccccccccc}
 L & \rightarrow & L {\otimes}_K L & \simeq &N_1 & \oplus & ... & \oplus & N_r   \\
\uparrow &  & \uparrow & & \uparrow & & & & \uparrow \\
k & \rightarrow & k[\Gamma] & \simeq & M_1 & \oplus & ... & \oplus & M_r
\end{array} \]
in which we recall that $N_1 = L$. We have in particular $ r = 6 $ in the last example because $L/K$ is a Galois extension for $S_3$. We check at once that $\Gamma (K') \subset K' $ when $K' = K (j) $ as we have $( e(j)=j, \sigma (j) = j, ...\tau \sigma (j) = j^2 = -(j+1))$ and thus ${\Gamma}' \lhd \Gamma$. On the contrary, we have $\Gamma (K') \notin K'$ when $K' = K(\eta)$ because $\sigma (\eta) = j \eta$ and thus ${\Gamma}' $ is not normal in $\Gamma$.  \\  \\   \\

\noindent
{\bf 3) ALGEBRAIC TOOLS}  \\

In all what follows, we shall consider unitaty rings $A, B, ... $ with elements $1, a, b, c, ...$ and full rings of quotients $K=Q(A), L=Q(B), ... $ which are fields when $A,B, ...$ are integral domains, that is do not have divisors of $0$, namely $a, b \in A, ab=0  \Rightarrow  a=0 $ or $b=0$. When $X$ is an algebraic set defined over a field $k$ of characteristic zero, that is $\mathbb{Q} \subset k$, we may introduce as usual the ring $A=k[X]$ of polynomial functions on $X$. In particular, if $G$ is an algebraic group defined over $k$, we may introduce the ring $R=k[G]$ as an algebra over $k$. We recall the basic axioms of $G$ that depend on the following objects:  \\
\[  \left\{ \begin{array}{cccccccl}
\epsilon & : & G \times G & \rightarrow & G & :  & (a,b) \rightarrow ab &  ({\it composition})  \\
\iota      & : &  G          & \rightarrow & G            &  : &       a \rightarrow a^{-1} & ({\it inverse}) \\
      e & : &  G & \rightarrow  & G  &  : &  a \rightarrow e & ({\it identity})
\end{array}  \right. \]
with the following axioms and corresponding commutative diagrams:  \\
\[     a, b, c \in G \rightarrow  (a b) c= a ( b c) = a b c \in G. \]
\[   \begin{array}{rcccl}
  & G \times G \times G & \stackrel{(\epsilon,{id}_G)}{\longrightarrow} & G \times G  & \\
({id}_G,\epsilon)           &  \downarrow  &          &       \downarrow    & \epsilon    \\
   &G \times G  &      \stackrel{\epsilon}{\longrightarrow }  &  G  &     
\end{array} \]  \\

\[  a \in G  \rightarrow  a^{-1} \in G , \hspace{1cm}  a a^{-1} = a^{- 1} a= e  \]
\[   \begin{array}{rcccl}
  &  G  & \stackrel{(\iota,{id}_G)}{\longrightarrow} & G \times G  & \\
({id}_G,\iota)           &  \downarrow  &  \searrow   e     &       \downarrow    & \epsilon    \\
   &G \times G  &      \stackrel{\epsilon}{\longrightarrow }  &  G  &     
\end{array} \]   \\

\[  e a = a e = a        \]   
\[   \begin{array}{rcccl}
  &  G  & \stackrel{( e, {id}_G)}{\longrightarrow} & G \times G  & \\
({id}_G, e)           &  \downarrow {id}_G &  \searrow  {id}_G     &       \downarrow    & \epsilon    \\
   &G \times G  &      \stackrel{\epsilon}{\longrightarrow }  &  G  &     
\end{array} \]   \\
where we have set $ G \stackrel{e}{\rightarrow } G: a \rightarrow e$ for the map sending any element $a \in G$ to $e \in G$.  \\

When $R$ is an algebra over $k$, we may consider the comorphisms of the morphisms appearing in the preceding diagrams and such a ring $R$ will be called a {\it Hopf algebra} over $k$. We may define:  \\

\[  \left\{ \begin{array}{ccccccl}
\tilde{\epsilon} & : & R & \rightarrow & R \, {\otimes}_k R & :  &  ({\it diagonal})  \\
\tilde{\iota}      & : &  R          & \rightarrow & R           &  : &        ({\it antipode}) \\
\tilde{e} & : &  R & \rightarrow  & k \subset R  &  : &   ({\it augmentation})
\end{array}  \right. \]

with the following commutative diagrams:  \\

\[   \begin{array}{rcccl}
  & R  & \stackrel{\tilde{\epsilon}}{\longrightarrow} & R \, {\otimes}_k R  & \\
 \tilde{\epsilon}           &  \downarrow  &          &       \downarrow    & {id}_R \otimes \tilde{\epsilon}    \\
   &R \, {\otimes}_k R  &      \stackrel{\tilde{\epsilon}\otimes {id}_R}{\longrightarrow }  & R \, {\otimes}_k R {\otimes}_k R &     
\end{array} \]  \\

\[   \begin{array}{rcccl}
  &  R  & \stackrel{\tilde{\epsilon}}{\longrightarrow} & R \, {\otimes}_k R  & \\
\tilde{\epsilon}          &  \downarrow  &  \searrow   \tilde{e}     &       \downarrow    & ({id}_R, \tilde{\iota})   \\
   &R \, {\otimes}_k R  &      \underset{(\tilde{\iota}, {id}_R)}{\longrightarrow }  &  R  &     
\end{array} \]   \\

\[   \begin{array}{rcccl}
  &  R  & \stackrel{\tilde{\epsilon}}{\longrightarrow} & R \, {\otimes}_k R  & \\
\tilde{\epsilon}          &  \downarrow  &  \searrow  {id}_R     &       \downarrow    & ({id}_R, \tilde{e})    \\
   &R \, {\otimes}_k R  &      \underset{(\tilde{e},{id}_R)}{\longrightarrow }  &  R  &     
\end{array} \]   \\
where we have set $(\tilde{\iota}, {id}_R): R \, {\otimes}_k R \stackrel{\tilde{\iota}{\otimes} {id}_R}{\longrightarrow} R \, {\otimes}_k R \rightarrow R$ and the last morphism is $a\otimes b \rightarrow ab$.  \\

These definitions are just those of the corresponding comorphisms when $R=k[G]$. However, the main difficulty in general is that one cannot find {\it any} generic solution of the defining finite Lie equations that could depend on certain arbitrary functions. The two examples of algebraic pseudogroups provided by the Pfaffian system ${\bar{y}}^2 d {\bar{y}}^1 = y^2 d y^1$ or by the Schwarzian OD equation $\frac{\partial \bar{y}}{\partial y} \frac{{\partial}^3 \bar{y}}{\partial y^3} - \frac{3}{2}(\frac{{\partial}^2 \bar{y}}{\partial y^2})^2=0$ are well known. The following key theorem provides the tricky differential geometric counterpart of the results of ([2, 3]) that will be essential in the next sections.  \\

\noindent
{\bf THEOREM 3.1}: If the manifold $X$ is a PHS for $G$, that is when the graph of the action $X \times G \rightarrow X\times Y$ is an isomorphism when $Y$ is a copy of $X$, the group parameters are {\it constants} on $X \times Y$ for the reciprocal distribution of the infinitesimal action of $G$ on $X$. Of course, the simplest example is that of a Lie group $G$ acting on itself. \\

\noindent
{\it Proof}: In the purely algebraic framework, counting the dimensions, we have $n=dim(X) = dim(G)=p$ and we may exhibit $p$ functions $a=\varphi(x,y)$ such that we have the $n$ identities $y\equiv f(x, \varphi(x,y))$. If $\delta = {\xi}^i(x)\frac{\partial}{\partial x^i} \in \Delta$ is a transformation commuting with all the infinitesimal generating transformations $\Theta = ({\theta}_1, ... , {\theta}_p)$, we can extend each $\delta$  to $X \times Y$ by setting $\delta = {\xi}^i(x) \frac{\partial}{\partial x^i} + {\xi}^k(y) \frac{\partial}{\partial y^k}$. Applying such an extended $\delta$ to the previous identities, we obtain:  \\
\[       {\xi}^k(y)= {\xi}^i(x)\frac{\partial f^k}{\partial x^i} (x,a) + (\delta {a}^{\tau} ) \frac{\partial f^k}{\partial a^{\tau}}(x,a) \] 
whenever $a=\varphi (x,y)$. As $[\delta, {\theta}_{\tau}]=0, \forall \tau = 1, ... ,p$, we have 
${\xi}(y)= {\xi }(x)\frac{\partial f}{\partial x}$. Also, as $X$ is a PHS for $G$, we ave $rk(\frac{\partial f}{\partial a})=n=p$ and thus 
$\delta a^{\tau} = 0$. Finally, when $\Phi$ is an invariant of $G$, that is $\theta \Phi=0$, as $[ \Delta , \Theta ]=0$, we get $ \theta (\delta \Phi)=\delta (\theta \Phi) =0$ and thus $\delta \Phi$ {\it must } be one invariant of $\Theta$.  \\
In the purely algebraic framework, let $X$ be an irreducible variety defined over the field $K$ and $G$ be an algebraic group defined over the field $k \subset K$. The group parameters must be therefore " {\it constants} " for $\Delta$, that is $ k[G] \subset cst( L {\otimes}_K L) \subset cst (Q(L {\otimes }_K L)) $. Accordingly, setting $A=K[X] \Rightarrow L = Q(A) = K(X)$ and introducing the subfield $K= L ^{G} \subset L$ invariant by $G$, we have thus an isomorphism $ Q(L {\otimes }_k k[G] )\simeq  Q( L {\otimes }_K L) $. With more details, $L {\otimes }_K L$ is a direct sum of integral domains acted on {\it separately} by $\Delta$ while  $Q(L {\otimes}_K L)$ is a direct sum of fields. In the case of a Lie group acting on itself, one has just to use the fact that {\it the left invariant distribution commutes with the right invariant distribution}.  \\
\hspace*{12cm}  $ \Box $      \\

Recapitulating these results while introducing the two injective comorphisms $L \rightarrow  L{\otimes}_K L$ respectively defined by 
$ \tilde{\alpha }: a \rightarrow a \otimes 1$ and $ \tilde{\beta}: a \rightarrow 1 \otimes a $ by analogy with the {\it source projection} $ \alpha :X \times Y \rightarrow X: (x,y) \rightarrow x$ and the {\it target projection} $ \beta : X \times Y \rightarrow Y : (x,y) \rightarrow y$, we discover that $L {\otimes}_K L$ is the simplest possible cogroupoid. We obtain indeed:   \\ \\
{\it Diagonal} : $k[G] \rightarrow k[G] {\otimes}_k k[G]$ is induced from : 
\[ L {\otimes}_K L \rightarrow (L {\otimes}_K L) {\otimes}_L ( L {\otimes}_K L ) \simeq   L {\otimes}_K L {\otimes}_K L : 
a \otimes b \rightarrow a \otimes 1 \otimes b \]
{\it Antipode} : $k[G] \rightarrow k[G]$ is induced by $  L {\otimes}_K L \rightarrow  L {\otimes}_K L : 
a \otimes b \rightarrow  b \otimes a $.  \\   \\
{\it Augmentation} : $k[G] \rightarrow  k \subset k[G]$ is induced by $ L {\otimes}_K L \rightarrow L : a \otimes b \rightarrow  a b $.  \\
if we identify $a\in L$ with $a \otimes 1$ and $b \in L$ with $1 \otimes b$ in $ L { \otimes}_K L$.  \\

\noindent
{\bf EXAMPLE 3.2}: With $n=1, m=2, k= \mathbb{Q}$, if $K$ is a differential field containing $k$, we may consider the differential automorphic extension $L/K$ defined by $L=K(y^1, y^2, y^1_x, y^2_x)$ in such a way that $d_x y^1=y^1_x, d_x y^2=y^2_x, d_x y^1_x=0, d_x y^2_x=0$. It amounts to consider the automorphic system $y^1_{xx}=0, y^2_{xx}=0$ under the standard action of $GL(2)$ which cannot be treated through the classical approach of the Picard-Vessiot theory introduced by Kolchin. In the present situation we have indeed ${\bar{y}}^1= a y^1 + b y^2, {\bar{y}}^2= c y^1 + d y^2$ and obtain easily the $4$ functions $\varphi$ with (See Example 4.26):
\[     \fbox{   $   a= (y^2_x {\bar{y}}^1 - y^2 {\bar{y}}^1_x)/( y^1 y^2_x - y^2 y^1_x)  $ }  \]
The $4$ infinitesimal generators of the action are $\Theta = \{ y^l \frac{\partial }{\partial y^k}\mid k,l=1,2 \}$ and their respective first order prolongations ([27]). Then $\Delta $ surely contains $\delta = y^k \frac{\partial}{\partial y^k_x}$ that we may extend to $\delta = y^k \frac{\partial}{ \partial y^k_x} + {\bar{y}}^k \frac{\partial}{\partial {\bar{y}}^k_x}$. We let the reader check that $\delta a = 0$, a result which is not evident at all. Of course, $L/K$ is {\it not} a PV extension in the sense of Kolchin because $y^1_x, y^2_x \in L$ are "{\it differential constants} " not in $K$. Accordingly, if we consider the intermediate differential field $K' = K(y^2_x)$ with the strict inclusions $K \subset K' \subset L$, the subgroup of invariance preserving $y^2_x$ is defined by ${\bar{y}}^2_x=c y^1 + d y^2$, that is we must have $c=0, d=1$ and this group also preserve the intermediate differential field $K'' = K(y^2, y^2_x)$ with the strict inclusion $K' \subset K'' $. It follows that we {\it cannot} have a Galois-type correspondence. Of course, such a result is coherent with the fact that $L/K'$ is NOT a PV differential extension in the sense of Kolchin because $L$ contains new differential constants as we have indeed $ y^1_x \in L, d_x y^1_x=0$. On the contrary, if we adopt the point of view of BB, we have $\delta y^2_x= y^2$ and $K'$ is not stable under $\Delta$ but it is a pity that the work of BB has {\it never} been acknowledged by Kolchin who had even never been able to use tensor products of rings and fields as any reader can check directly through {\it all} his books.  \\   \\

\noindent
{\bf 4) DIFFERENTIAL TOOLS}:  \\

Our purpose is now to use the fact that Lie groups of transformations are just exampes of Lie pseudogroups of transformations in such a way that $trd(L/K)=\infty$. However, differentiating the group law $y=f(x,a)$ in order to obtain the prolongations $y_x= {\partial}_x f(x,a), y_{xx}= {\partial}_{xx} f(x,a)$ and so on, eliminating the parameters may be quite a hard task like in the case of the projective group of transformations of the real line (Exercise). Also, no classical method known for Lie groups can work for Lie pseudogroups because one cannot find generic solutions in most cases when $ trd(L/K)=\infty $.   \\

We start with a few basic technical results and formulas that are not well known because they involve the Spencer operator. For simplicity, we shall deal with trivial fibered manifolds ${\cal{E}}=X \times Y$ such that $dim(X)=n, dim(Y)=m$ and local coordinates $ (x^i, y^k) $ with $ i=1,...,n, k=1,...,m $. The $q$-jet bundle $J_q(\cal{E}) $ of $\cal{E}$ will be a fibered manifold with local coordinates $(x^i, y^k_{\mu})$ for a multi-index $\mu = ({\mu}_1, ..., {\mu}_n)$ of length $0\leq  \mid \mu \mid={\mu}_1+ ...+ {\mu}_n \leq q$ and we shall set $ \mu + 1_i= ({\mu}_1, ...,{ \mu}_{i-1} , {\mu}_i+1,{\mu}_{i +1 }, ..., {\mu}_n) $ or simply $(x,y_q) $ with projection $\pi$ to $X$.  The tangent bundle $T(\cal{E}) $ may be described by means of local coordinates $ (x,y; u,v)$ while the vertical bundle $V(\cal{E}) $ will be obtained by setting $u=0$ in the short exact sequence of vector bundles pulled back over 
$\cal{E}$:  \\
\[ 0   \rightarrow  V({\cal{E}}) \rightarrow T({\cal{E}}) \stackrel { T(\pi) }{ \longrightarrow} T {\times}_X {\cal{E}} \rightarrow 0 \]
Introducing the formal derivatives $d_i= {\partial}_i + y^k_{\mu + 1_i}\frac{\partial }{\partial y^k_{\mu}} $, we have (See [27-30, 35] for details): \\

\noindent
{\bf LEMMA 4.1}: Prolongation of vertical vector fields:  \\
\[\eta= {\eta}^k(y) \frac{\partial}{\partial y^k}   \rightarrow  {\rho}_q(\eta)= d_{\mu}{\eta}^k\frac{\partial}{\partial y^k_{\mu}} \]
\[ d_i{\eta}^k(y)= \frac{\partial {\eta}^k}{\partial y^r} y^r_i, d_{ij}{\eta}^k(y)= \frac{\partial {\eta}^k}{\partial y^r}y^r_{ij} + \frac{{\partial}^2 {\eta}^k}{\partial y^r \partial y^s}y^r_i y^s_j , ....  \]
Introducing a section ${\eta}_q$ of $J_q(T(Y))$ over the target and replacing derivatives by sections, we get:  \\
\[ \fbox{  $  \sharp({\eta}_q)= {\eta}^k\frac{\partial}{\partial y^k} + {\eta}^k_r y^r_i \frac{\partial }{\partial y^k_i} + ({\eta}^k_r y^r_{ij} + {\eta}^k_{rs}y^r_iy^s_j) \frac{\partial }{\partial y^k_{ij}} + ... $  }  \]

\noindent
{\bf LEMMA 4.2}: Prolongation of horizontal vector fields:  \\
\[\xi= {\xi}^i(x) {\partial}_i  \rightarrow  {\rho}_q(\xi)= {\xi}^i {\partial}_i + {\zeta}^k_{\mu} \frac{\partial}{\partial y^k_{\mu}} \]
\[ {\zeta}^k_{\mu + 1_i} = d_i{\zeta}^k_{\mu} - y^k_{\mu +1_r}{\partial}_i{\xi}^r(x) \]
\[ {\zeta}^k=0, \,\, {\zeta}^k_i=  - y^k_r {\partial}_i {\xi}^r, \,\, 
{\zeta}^k_{ij}=  - (y^k_r {\partial}_{ij} {\xi}^r + y^k_{r j}{\partial}_i{\xi}^r + y^k_{r i} {\partial}_j{\xi}^r) , ....  \]
Introducing a section ${\xi}_q$ of $J_q(T)$ over the source and replacing derivations by sections, we get:  \\
\[ \fbox{  $  \flat({\xi}_q)= {\xi}^i(x){\partial}_i - y^k_r {\xi}^r_i(x) \frac{\partial }{\partial y^k_i} - (y^k_r {\xi}^r_{ij}(x) + 
y^k_{r j} {\xi}^r _i(x)+ y^k_{r i} {\xi}^r_j(x)) \frac{\partial }{\partial y^k_{ij}} + ... $  }   \]

\noindent
{\bf THEOREM 4.3}: There exists a bracket for sections of $J_q(T)$ generalizing the standard bracket of vector fields of $T$. \\

\noindent
{\it Proof}: We recall that $([\xi,\eta])^i={\xi}^r{\partial}_r{\eta}^i - {\eta}^r{\partial}_r{\xi}^i, \forall \xi, \eta \in T$. Taking the $q$-derivation by applying the operator $j_q$, we obtain a bilinear combination of $j_{q+1}(\xi)$ and $j_{q+1}(\eta)$. We may thus define the so-called {\it algebraic bracket} $\{ {\xi}_{q+1}, {\eta}_{q+1}\}$ with value in $J_q(T)$ and obtain on $J_q(T)$ the {\it algebroid bracket} where $d $ is the Spencer operator $ (d {\xi}_{q+1})^k_{ \mu, i} (x)= {\partial}_i {\xi}^k_{\mu} (x)- {\xi}^k_{ \mu + 1_i}(x)$ or simply $d{\xi}_{q+1} = j_1({\xi}_q ) - {\xi}_{q+1}$: \\
\[  \fbox{   $ [{\xi}_q, {\eta}_q ] = \{ {\xi}_{q+1}, {\eta}_{q+1} \} + i(\xi) d {\eta}_{q+1} - i(\eta) d {\xi}_{q+1}  $  }   \]
which does not depend any longer on the jets of strict order $q+1$ whenever ${\xi}_{q+1}, {\eta}_{q+1} \in J_{q+1}(T)$ are over 
${\xi}_q, {\eta}_q \in J_q(T)$. When ${\xi}_{q+1}=j_{q+1}(\xi)$ and ${\eta}_{q+1}=j_{q+1}(\eta)$, we have $d{\xi}_{q+1}=0, d{\eta}_{q+1}=0$ 
and thus $[j_q(\xi), j_q(\eta)]=j_q([\xi,\eta]) $ .
It is finally highly not evident to verify the Jacobi identity:  \\
\[  \fbox{  $  [ {\xi}_q , [{\eta}_q,{\zeta}_q] ] + [ {\eta}_q, [{\zeta}_q, {\xi}_q] ] +  [ {\zeta}_q, [ {\xi}_q, {\eta}_q] ] =0 , 
\forall {\xi}_q, {\eta}_q, {\zeta}_q \in J_q(T) $  }  \]
\hspace*{12cm}     $ \Box  $   \\

\noindent
{\bf DEFINITION 4.4}: A sub-bundle $R_q \subset J_q(T)$ is a {\it Lie algebroid} of order $q$ if $[R_q,R_q] \subset R_q$, that is $[{\xi}_q, {\eta}_q] \subset R_q , \forall {\xi}_q, {\eta}_q \in R_q$. We say that $R_q$ is {\it transitive} if the morphism $R_q \rightarrow T$ induced by the canonical epimorphism ${\pi}^q_0: J_q(T) \rightarrow T$ is also an epimorphism. We shall introduce the {\it isotropy Lie algebra bundle} $R^0_q$ by the short exact sequence $0 \rightarrow R^0_q \rightarrow  R_q \rightarrow T \rightarrow 0$. We have $[R^0_q, R^0_q] \subset R^0_q$ {\it fiber by fiber} and a $R_q$-connection is a map ${\chi}_q:T \rightarrow R_q$ such that ${\pi}^q_0 \circ {\chi}_q = id_T$.\\

\noindent
{\bf COROLLARY 4.5}: One has $[ \sharp({\eta}_q), \sharp({\eta}'_q] = \sharp ([{\eta}_q,{\eta}'_q])$ over the target and 
$[\flat ({\xi}_q), \flat({\xi}'_q)]=\flat ([{\xi}_q,{\xi}'_q])$ over the source.  \\

\noindent
{\bf COROLLARY 4.6}: One has $[\flat({\xi}_q), \sharp({\eta}_q)]=0, \forall {\xi}_q \in J_q(T)$ over the source, $\forall{\eta}_q \in J_q(T(Y))$ over the target as a generalization for sections of the fact that any source transformation commutes with any target transformation, namely that $[\xi,\eta]=0, \forall \xi \in T=T(X), \forall \eta \in T(Y)$.  \\

In actual practice, apart from ([21, 34]), we do not know any reference on these results which are crucially depending on the use of the Spencer operator. We invite threader to check that a direct proof of these formulas is rather easy when $q=1$ but becomes quite tricky even when $q=2$. The same comment is valid for the following two formulas which are among the most difficult but also the most useful ones. \\

\noindent
{\bf THEOREM 4.7}: For any function $\Phi(x,y_q)$, we have:  \\
\[ \fbox{  $   \sharp({\eta}_{q+1}) d_i \Phi= d_i (\sharp({\eta}_q) \Phi) - y^k_i (d_Y{\eta}_{q+1} (\frac{\partial}{\partial y^k})) \Phi  $ }  \]
when ${\eta}_{q+1}$ projects onto ${\eta}_q$ over the target and the similar formula over the source:  \\
\[ \fbox{ $  \flat({\xi}_{q+1} ) d_i \Phi = d_i( \flat ({\xi}_q) \Phi) - {\xi}^r_i d_r \Phi  - \flat(d {\xi}_{q+1}({\partial}_i))\phi   $  }  \]

\noindent
{\bf COROLLARY 4.8}: If $ \Phi=\Phi(y_q)$ is a differential invariant of strict order $q$ and ${\eta}_{q+1}\in R_{q+1}(Y)$ over the target, then $d_Y{\eta}_{q+1} \in T^*(Y) \otimes R_q(Y)$ over the target and $d_i \Phi$ is thus a differential invariant of strict order $q+1$.  \\

\noindent
{\bf EXAMPLE 4.9}: In order to help the reader dealing with sections instead of solutions, we consider the case $n=m=1$ with 
${\xi}_2=(\xi, {\xi}_x, {\xi}_{xx})$ and we have:
\[ \flat({\xi}_2)= \xi {\partial}_x + 0 \frac{\partial}{\partial y} - y_x {\xi}_x \frac{\partial}{\partial y_x} - 
(y_x {\xi}_{xx} + 2 y_{xx}{\xi}_x) \frac{\partial}{\partial y_{xx}}   \]
We let the reader prove, as a tricky exercise, that:  \\
\[   \flat({\xi}_2) d_x \Phi = d_x (\flat({\xi}_1) \phi) - {\xi}_x d_x \phi - \flat (d{\xi}_2({\partial}_x))\Phi \]
When ${\xi}_2= (0,-1,0) \Rightarrow d{\xi}_2 =(1,0)$ and $\Phi=\Phi (y,y_x)$, then one has:  \\
\[    (y_x \frac{\partial}{\partial y_x} + 2 y_{xx} \frac{\partial}{\partial y_{xx}}) (y_x \frac{\partial \Phi}{ \partial y} + y_{xx} \frac{\partial \Phi}{\partial y_x})=
       d_x (y_x \frac{\partial \Phi}{\partial y_x})  + d_x \Phi \]

\noindent
{\bf REMARK 4.10}: If $\eta$ is an infinitesimal generator of the action of a Lie group on $Y$, choosing ${\eta}_{q+1} = j_{q+1}(\eta)$ that projects onto $  {\eta}_q=j_q(\eta)$ over the target, we obtain $ {\rho}_{q+1}(\eta) d_i \Phi=d_i ({\rho}_q (\eta) \Phi)$. It follows that $d_i$ commutes with the prolongations of target transformations as we already checked on examples.  \\

\noindent
{\bf THEOREM 4.11}: We have the following formula for the bracket of sections of $ J_{q+1}(T)$:  \\
\[   i(\zeta) d[{\xi}_{q+1}, {\eta}_{q+1}] =  [ i(\zeta) d{\xi}_{q+1}, {\eta}_{q+1}] + [{\xi}_q, i(\zeta) d{\eta}_{q+1}]  \]
for {\it any} $\zeta \in T$ (See [27, 35] for more details).  \\

By an induction starting with $r=1$, we obtain:  \\

\noindent
{\bf COROLLARY 4.12}: If $[R_q,R_q] \subset R_q$, then $[R_{q+r},R_{q+r}] \subset R_{q+r} , \forall r \geq 0$ and the prolongations of a Lie algebroid are Lie algebroids even if $R_q \subset J_q(T)$ is not formally integrable. The case of the Killing system for the Schwarzschild and Kerr metrics provides a good example in general relativity.  \\

\noindent
{\bf EXAMPLE 4.13}: When $m=n=2$, the algebraic pseudogroup defined by the Pfaffian system $y^2 dy^1=x^2 dx^1$ is equivalently defined by the nonlinear system $y^2 \frac{\partial y^1}{\partial x^1} =x^2, \frac{\partial y^1}{\partial x^2}=0$. The corresponding Lie algebroid $R_1 \subset J_1(T)$ is defined by the linearized system $ x^2 {\xi}^1_1 + {\xi}^2=0, {\xi}^1_2=0$ which is not involutive and not even formally integrable because, using crossed derivatives, we obtain at once ${\xi}^1_1 + {\xi}^2_2=0$. Hence, we can use two sections ${\xi}_1, {\eta}_1 \in R_1$ with ${\xi}^1_1 +{\xi}^2_2\neq 0, {\eta}^1_1 + {\eta}^2_2 \neq 0$ and check that $[{\xi}_1,{\eta}_1]\in R_1$ in such a way that $({\xi}^1=0, {\xi}^2=-x^2, {\xi}^1_1=1, {\xi}^1_2=0, {\xi}^2_1=0, {\xi}^2_2=0) $ and $({\eta}^1=1, {\eta}^2=-0, {\eta}^1_1=0, {\eta}^1_2=0, {\eta}^2_1=1, {\eta}^2_2=1) $ {\it with no relation at all with solutions}. \\

Considering now an algebraic pseudogroup defined {\it over the target} by a nonlinear system ${\cal{R}}_q(Y) \subset {\Pi}_q(Y,Y) \subset J_q(Y\times Y)$ with local coordinates simply denoted by $(y, \bar{y}, \frac{\partial \bar{y}}{\partial y}, \frac{{\partial}^2 \bar{y}}{\partial y \partial y}, ...) $ such that $det( \frac{\partial \bar{y}}{\partial y}) \neq 0$. The jet composition $((\bar{y},\bar{\bar{y}},\frac{\partial \bar{\bar{y}}}{\partial \bar{y}}, ...)(y,\bar{y},\frac{\partial \bar{y}}{\partial y}, ...)) \rightarrow   (y, \bar{\bar{y}}, \frac{\partial \bar{\bar{y}}}{\partial \bar{y}}\frac{\partial \bar{y}}{\partial y}, ...)$ is obtained by using the chain rule for derivatives while the inversion is $ (y, \bar{y}, \frac{\partial \bar{y}}{\partial y}, ... ) \rightarrow ( \bar{y}, y, (\frac{\partial \bar{y}}{\partial y})^{-1}, ...) $.  \\
The composition $J_q(X\times Y) {\times}_Y{\cal{R}}_q(Y) \rightarrow J_q(X\times Y)$ may be similarly defined by using:  \\
\[   (x,y,y_x, y_{xx}, ...) (y, \bar{y}, \frac{\partial \bar{y}}{\partial y}, \frac{{\partial}^2 \bar{y}}{\partial y \partial y}, ...)  \rightarrow
 (x, \bar{y}, \frac{\partial \bar{y}}{\partial y} y_x, \frac{\partial \bar{y}}{\partial y} y_{xx} +
  \frac{{\partial}^2 \bar{y}}{\partial y \partial y} y_x y_x, ...)\]
  By a {\it free generic action} we understand that the morphism $\sharp: R_q(Y) \rightarrow J_q(V(X \times Y))\simeq V(J_q(X \times Y))$ is a monomorphism and we can thus set:  \\

\noindent
{\bf DEFINITION 4.14}: A nonlinear system ${\cal{A}}_q \subset J_q(X\times Y) $ over the source is said to be a PHS for the Lie groupoid ${\cal{R}}_q(Y) \subset {\Pi}_q(Y,Y)$ over the target if the corresponding graph:  \\
\[ {\cal{A}}_q {\times}_Y {\cal{R}}_q(Y) \rightarrow {\cal{A}}_q \times {\cal{A}}_q    \]
is an isomorphism. Setting $A_q=V({\cal{A}}_q)$, the action is said to be free (transitive, simply transitive) when the morphism $\sharp : {\cal{A}}_q {\times}_Y R_q(Y) \rightarrow A_q$ of vector bundles over ${\cal{A}}_q$ is a monomorphism (an epimorphism, an isomorphism). The system is said to be an {\it automorphic system} if the $r$-prolongation ${\cal{A}}_{q+r} $ is a PHS for the $r$-prolongation ${\cal{R}}_{q+r}(Y)$ , $\forall r\geq 0$.  \\

We have already proved and illustrated in ([28, 36]) the two following criteria:  \\

\noindent
{\bf THEOREM 4.15}: ({\it First criterion for automorphic systems})  If an involutive system $ {\cal{A}}_q \subset J_q(X\times Y)$ is a PHS for a Lie groupoid ${\cal{R}}_q(Y) \subset {\Pi}_q(Y,Y) $ and if $ {\cal{A}}_{q+1} ={\rho}_1({\cal{A}}_q) \subset J_{q+1}(X\times Y)$ is a PHS for the Lie groupoid ${\cal{R}}_{q+1}(Y) ={\rho}_1({\cal{R}}_q(Y)) \subset {\Pi}_{q+1}(Y,Y) $, then ${\cal{R}}_q$ is an involutive system over the target with the same non-zero characters and ${\cal{A}}_q$ is an automorphic system.  \\

\noindent
{\bf THEOREM 4.16}: ({\it Second criterion for automorphic systems}) If ${\cal{R}}_q \subset {\Pi}_q(Y,Y)$ is an involutive system of finite Lie equations such that the action of ${\cal{R}}_q(Y)$ on $J_q(X\times Y)$ is generically free, then the action of ${\cal{R}}_{q+r}(Y) $ on $J_{q+r}(X\times Y)$ is generically free $\forall r \geq 0$ and all the differential invariants are generated by by a fundamental set of order $q+1$ ({\it care}).  \\

\noindent
{\bf EXAMPLE 4.17}: when $n=1,m=2, q=2$ and $k=\mathbb{Q}$, we may consider the Lie pseudogroup $\Gamma$ defined as a Lie group by the action $\bar{y}=Ay + B$ with $det(A)=1$. The {\it only} generating differential invariant at order $2$ is $\Phi= y^1_xy^2_{xx} - y^2_xy^1_{xx}$ but, at order $3$, we must use $d_x\Phi=y^1_x y^2_{xxx} - y^2x y^1_{xxx} $ of course but we have also to add $ \Psi = y^1_{xx} y^2_{xxx} - y^2_{xx} y^1_{xxx}$. We have thus the {\it strict inclusion} ${\cal{A}}_3 \subset {\rho}_1({\cal{A}}_2)$ and the symbol of ${\cal{A}}_3$ is vanishing if and only if $y^1_x y^2_{xx} - y^2_x y^1_{xx} \neq 0$. We shall meet a similar condition with the non-zero Wronskian determinant $y^1 y^2_x  -  y^2 y^1_x\neq 0 $ at order $q=1$ in the Picard-Vessiot theory if we consider the action $\bar{y}=Ay$ with $det(A)=1$. Indeed, we must consider the only first order differential invariant 
$ \Phi= y^1 y^2_x -y^2 y^1_x$ and use the second order differential invariant $d_x\Phi= y^1 y^2_{xx} - y^2 y^1_{xx}$ to which we must add $\Psi= y^1_x y^2_{xx} - y^2_x y^1_{xx}$. We notice that the symbol of order $2$ is vanishing if and only if $y^1 y^2_x - y^2 y^1_x \neq 0 $ and let the reader check that all these results are coherent with the two previous criteria.  \\

We are now able to recognize whether a nonlinear system of algebraic OD or PD equations with $n$ independent variables and $m$ unknowns is an automorphic system for its biggest pseudogroup of invariance. It thus remains to exhibit the transition to differential algebra. For this, keeping in mind the difference existing between special and general relativity in physics, we shall explain on an example the difference existing between a special and a general automorphic system.  \\

\noindent
{\bf EXAMPLE 4.18}: With $n=1, m=2, q=1, k=\mathbb{Q}$, let us consider the algebraic pseudogroup of target transformations preserving the $1$-form $y^2dy^1$ and thus also the $2$-form $d y^1 \wedge d y^2$. On one side we may start with a given differential field $K$ containing $k$ and consider the special system $y^2 y^1_x= \omega \in K$. Looking for the biggest lie pseudogroup of invariance of this OD equation, we must have:  \\
\[   {\bar{y}}^2 {\bar{y}}^1_x={\bar{y}}^2 \frac{\partial {\bar{y}}^1}{ \partial y^1} y^1_x + {\bar{y}}^2 \frac{\partial {\bar{y}}^1}{ \partial y^2} y^2_x = y^2 y^1_x \Rightarrow {\bar{y}}^2 \frac{\partial {\bar{y}}^1}{\partial y^1} = y^2, \frac{\partial {\bar{y}}^1}{\partial y^2}=0 \]
Such a system is not involutive as it is not even formally integrable and we must add 
$\frac{\partial ({\bar{y}}^1,{\bar{y}}^2)}{\partial (y^1,y^2)}=1$. We may thus start with the linear involutive system $R_1(Y)$ defined by the first order system of infinitesimal Lie equations $y^2 {\eta}^1_1 + {\eta}^2=0, {\eta}^1_2=0, {\eta}^1_1 + {\eta}^2_2=0$. Introducing the prime differential ideal $\mathfrak{p} \subset K \{y\}=K[y,y_x,y_{xx}, ...] $ we may define the {\it special} differential automorphic extension $L/K$ with $L= Q(K\{y\}/ \mathfrak{p}) \simeq K(y^1,y^2, y^2_x, y^2_{xx}, ...)$. However, we may also start with the differential extension $L=Q(k\{y\})=k<y>$ and consider the {\it general} differential automorphic extension $L/K$ with 
$ K=k<\Phi > $ by introducing the generating differential invariant $\Phi = y^2 y^1_x$.  \\

With a slight abuse of language, we may set ${\Theta}(q)=\sharp (R_q(Y))$ and this distribution can be studied without introducing the space of solutions of the Lie algebroid $R_q(Y)$ that may not be known explicitly, for example if one replaces $ y^2d y^1 $ by $y^2 d y^1 - y^1 d y^2$. Using the Frobenius theorem for the involutive distribution ${\Theta}(q)$ on $V(J_q(X \times Y)$, we may find a {\it reciprocal distribution} ${\Delta}(q)$ of maximal rank such that $[ {\Theta}(q), {\Delta}(q) ] =0$, that is $[\theta, \delta ]=0, \forall \theta \in \Theta, \delta \in \Delta$.  \\

\noindent
{\bf LEMMA 4.19}: $ \delta y = 0, \forall \delta \in  \Delta $.  \\

\noindent
{\it Proof}: As the sections of $R_q(Y)$ are only defined up to a linear combination with coefficients only depending on $y$, we have successively:  \\
\[  {\theta}' = \Sigma \lambda (y) \theta \Rightarrow  [\delta , {\theta}' ]= \Sigma (\delta \lambda) \theta + \Sigma \lambda [\delta , \theta ]= \Sigma (\delta \lambda) \theta = 0 \Rightarrow \delta \lambda=0  \]
Quantities killed by $\Delta$ will be called " {\it constants} " according to ([2, 3]) and we have thus $cst(L)=k(y)$. The next example will illustrate this new concept already 
introduced ... {\it forty years ago} in ([28]) but never acknowledged ([20, 21]). These results bring the need to revisit entirely the Picard-Vessiot theory as it is known today.  \\
\hspace*{12cm} $\Box $ \\

\noindent
{\bf EXAMPLE 4.20}: Coming back to the previous example, we get with $q=1$ and parametric sections $\{ {\eta}^1, {\eta}^2, {\eta}^2_1 \}$ of $R_1(Y)$, we get:  \\
\[  \fbox{  $   {\Theta}(1) \left \{ \begin{array}{ccccl}
{\eta}^2_1 & \rightarrow & {\theta}_3 & = &  y^1_x \frac{\partial }{ \partial y^2_x} \\
{\eta}^2     & \rightarrow & {\theta}_2 & = & y^2 \frac{\partial}{\partial y^2 } - y^1_x \frac{\partial}{\partial y^1_x} + y^2_x \frac{\partial}{\partial y^2_x} \\
{\eta}^1     & \rightarrow & {\theta}_1 & = &  \frac{\partial}{\partial y^1} 
\end{array} \right.    $  }   \]
We notice that the relation $y^1 {\theta}^1 - {\theta}^2= {\rho}_1 (y^1 \frac{\partial}{\partial y^1} - y^2 \frac{\partial}{\partial y^2}) $ proves that the pseudogroup of invariance contains finite transformations of the form $({\bar{y}}^1= a y^1, {\bar{y}}^2= \frac{1}{a}y^2)$ whenever $ a=cst$.  \\
A reciprocal distribution could be {\it a priori}:  \\
\[   \left \{ \begin{array}{ccc}
{\delta}_1     & = &  y^1_x \frac{\partial }{ \partial y^1_x} + y^2_x \frac{\partial }{ \partial y^2_x}\\
{\delta}_2     & = & {\hspace{8mm}}y^2 \frac{\partial}{\partial y^2_x}  \\
{\delta}_3     & = &  y^2  \frac{\partial}{\partial y^2} + y^2_x \frac{\partial}{\partial y^2_x} \\
{\delta}_4     & = &  \frac{\partial}{\partial y^1} { \hspace{13mm} }
\end{array} \right. \]
However, {\it taking into account the previous lemma}, we must push out ${\delta}_3 $ and ${\delta}_4$ in order to get:  \\
\[  \fbox{  $   {\Delta}(1) \left \{ \begin{array}{ccr}
{\delta}_1      & = &  y^1_x \frac{\partial }{ \partial y^1_x} + y^2_x \frac{\partial }{\partial y^2_x}  \\
{\delta}_2      & = &  y^2 \frac{\partial}{\partial y^2_x}  
\end{array} \right.   $   }  \]
which is of rank $2$ if and only if $y^2 y^1_x \neq 0$. We notice that ${\delta}^1$ is the factor of  $ - {\xi}_x$ in $\flat ({\xi}_1)$.\\
Similarly but after a tricky computation ({\it care}), we should obtain at order $2$:  \\
\[  \fbox{  $  {\Theta}(2) \left \{ \begin{array}{ccccc}
{\eta}^2_{11} & \rightarrow & {\theta}_4 & = & (y^1_x)^2 \frac{\partial}{\partial y^2_{xx}} \\
{\eta}^2_1     & \rightarrow & {\theta}_3 & = & y^2 y^1_x \frac{\partial }{ \partial y^2_x} - (y^1_x)^2 \frac{\partial}{\partial y^1_{xx}} +
 (y^2 y^1_{xx} + 2 y^1_x y^2_x)\frac{\partial}{\partial y^2_{xx}}\\
{\eta}^2         & \rightarrow & {\theta}_2 & = & y^2 \frac{\partial}{\partial y^2 } - y^1_x \frac{\partial}{\partial y^1_x} + y^2_x \frac{\partial}{\partial y^2_x} - y^1_{xx} \frac{\partial}{\partial y^1_{xx}} + y^2_{xx} \frac{\partial}{\partial y^2_{xx}} \\
{\eta}^1         & \rightarrow & {\theta}_1 & = &  \frac{\partial}{\partial y^1} 
\end{array} \right.    $  }  \]
and we have now $y^1 {\theta}_1 - {\theta}_2= {\rho}_2(y^1 \frac{\partial}{ \partial y^1} - y^2 \frac{\partial}{\partial y^2} )$.  \\
A reciprocal distribution could be:  \\
\[ \fbox{  $  {\Delta}(2)  \left \{ \begin{array}{ccc}
{\delta}_1     & = &  y^1_x \frac{\partial }{ \partial y^1_x} + y^2_x \frac{\partial }{ \partial y^2_x} + 2(y^1_{xx} \frac{\partial}{\partial y^1_{xx}} + y^2_{xx}  \frac{\partial}{\partial y^2_{xx}} )  \\
{\delta}_2     & = & y^2 \frac{\partial}{\partial y^2_x } + 2 y^2_x \frac{\partial}{\partial y^2_{xx}}   \\
{\delta}_3     & = &  y^2  \frac{\partial}{\partial y^2_{xx}}  \\
{\delta}_4     & = &  y^1_x \frac{\partial}{\partial y^1_{xx}} + y^2_x \frac{\partial}{\partial y^2_{xx}}
\end{array} \right.  $  }   \]
It is of rank $4$ in $(y^1_x, y^2_x, y^1_{xx}, y^2_{xx})$ and ${\delta}_4$ is now the factor of $ - {\xi}_{xx}$ in $\flat({\xi}_2) $.  \\
One can finally easily check that ${\Delta}(2)$ stabilizes both $\Phi=y^2 y^1_x$ and $d_x\Phi= y^2 y^1_{xx} + y^1_x y^2_x$.  \\

Passing from groups to pseudogroups while using the results of ([2, 3]), we have:  \\

\noindent
{\bf PROPOSITION 4.21}: Only stable intermediate differential fields can provide a Galois differential correspondence.  \\

\noindent
{\it Proof}: Let us consider the chain of inclusions $k \subset K \subset K' \subset L$ where the differential fields involved are respectively $ K=k< \Phi > $, $K'= k< \Phi, \Psi >$, $L= k< y >$ for simplicity. The pseudogroup ${\Gamma}' $ preserving $\Phi$ {\it and} $\Psi$ is thus a sub-pseudogroup of the pseudogroup $\Gamma$ preserving only $\Phi$. Hence, any new infinitesimal target transformation is of the form ${\theta}'= \Sigma {\lambda }(y) \theta$ with ${\theta}' \Psi =0$. We obtain therefore from the previous lemma $[ {\theta}' , \delta ] = \Sigma {\lambda }(y) [ \theta , \delta ] = O $ and thus $ {\theta}' (\delta \Psi)= \delta ({\theta}' \Psi) =0 $, that is to say $ \delta \Psi$ is {\it necessarily}  a differential invariant of ${\Gamma}' $, thus a differential function of $\Phi$ {\it and } $\Psi$.  \\
\hspace*{12cm} $\Box $ \\

\noindent
{\bf EXAMPLE 4.22}: If one chooses $k=\mathbb{Q}$ and $ K' = k<\Phi, \Psi >$ like in the previous example, that is with 
$\Phi= y^2 y^1_x$ and $ \Psi = y^2_x$, we obtain ${\delta}_1 \Psi= \Psi $ {\it but} ${\delta}_2 \Psi= y^2$, that is $K' $ is not stable. Indeed, the pseudogroup preserving $\Phi $ and $ \Psi$ is thus also preserving $y^2$, that is becomes ${\bar{y}}^1=y^1 + b, {\bar{y}}^2=y^2$. On the contrary, if we choose $\Psi= \frac{y^2_x}{y^2}$, then ${\delta}_1\Psi= =\Psi, , {\delta}_2 \Psi= =1$ and ${\Gamma}'$ becomes ${\bar{y}}^1= ay^1 + b, {\bar{y}}^2= \frac{1}{a} y^2$.  \\

{\it Last but not least}, we study the condition for which not only we shall have ${\Gamma}' \subset \Gamma$ but also ${\Gamma}' \lhd \Gamma$. First of all, up to the knowledge of the author, there does not exist formal conditions. Indeed, if $g \in \Gamma, h \in {\Gamma}'$, then $g^{-1} \circ h  \circ g = h' \in {\Gamma}'${\it cannot} be checked in general. \\

Now, the vector bundle $J_q(T)$ is also a bundle of geometric objects associated with ${\Pi}_{q+1}(X,X)$ in such a way that 
${\eta}_q=f_{q+1}({\xi}_q) \in J_q(T), \forall {\xi}_q \in J_q(T), \forall f_{q+1} \in {\Pi}_{q+1}$ along the formulas deduced from the transformations of a vector field and its various derivatives, namely ([NLCEM]):  \\
\[ {\eta}^k(f(x))= f^k_r(x){\xi}^r(x), {\eta}^k_u f^u_i = f^k_r {\xi}^r_i + f^k_{ri}  {\xi}^r ,     
   {\eta}^k_{uv}f^u_if^v_j+{\eta}^k_uf^u_{ij}  =  f^k_r{\xi}^r_{ij} + f^k_{ri}{\xi}^r_j+f^k_{rj}{\xi}^r_i+f^k_{rij}{\xi}^r \]
and so on. We may define the {\it finite normalizer} ${\tilde{\cal{R}}}_{q+1} = \{ {\tilde{f}}_{q+1} \mid {\tilde{f}}_{q+1} ({\xi}_q) \in R_q(Y) , \forall {\xi}_q \in R_q =R_q(X)  \}$ when $Y$ is a copy of $X$ {\it in this framework}.  \\

Passing to the infinitesimal point of view, we may define the {\it formal Lie derivative}:  \\
\[ L({\xi}_{q+1}){\eta}_q= [{\xi}_q, {\eta}_q ] + i(\eta) d {\xi}_{q+1} = \{ {\xi}_{q+1}, {\eta}_{q+1} \} + i(\xi) d {\eta}_{q+1}  \]
and define the {\it infinitesimal normalizer} by the condition ${\tilde{R}}_{q+1}= \{ {\tilde{\xi}}_{q+1} \mid L({\tilde{\xi}}_{q+1}) {\xi}_q \in R_q, \forall {\xi}_q \in R_q \}$
We obtain in particular the following quite difficult result ([27], p 390-393):  \\

\noindent
{\bf THEOREM 4.23}: If $R_q \subset J_q(T)$ is a formally integrable and transitive system of infinitesimal Lie equations with a symbol 
$g_q= R_q \cap S_qT^* \otimes T$ which is $2$-acyclic (involutive), then $ {\tilde{R}}_{q+1}$ is formally integrable (involutive) and such that ${\tilde{g}}_{q+1}=g_{q+1}$.  \\ 

\noindent
{\bf EXAMPLE 4.24}: Again with $n=1, m=2, k= \mathbb{Q}$, let us consider the inclusions $ k \subset K=k<\Psi> \subset K'=k<\Phi> \subset k<y^1, y^2 >$ with $ \Phi = y^2 y^1_x$ and $\Psi= (y^2 y y^1_{xx} + y^1_x y^2_x)/y^2 y^1_x= d_x \Phi / \Phi$. By shear luck, we can exhibit ${\Gamma}' = \{ {\bar{y}}^1=g(y^1), {\bar{y}}^2=  y^2 / (\frac{\partial g}{ \partial y^1} )\} \lhd \Gamma = \{ {\bar{y}}^1 = g(y^1), {\bar{y}}^2= a y^2 / (\frac{ \partial g}{\partial y^1}) \} $. According to the previous results, a necessary condition for having $\Gamma \lhd {\Gamma}' $ is that $diff trd (K' / K)=0$ because $diff trd (L/K) = diff trd (L/K')=1$ and $K'/K $ is a PV-extension for the multiplicative group of the real line.  \\

In view of the preceding results, even the general Picard-Vessiot theory must be entirely revisited in a coherent way with BB as follows.  \\

\noindent 
{\bf REMARK 4.25}:  The classical textbook definition $ N(\Theta) = \{ \xi \in T \mid [ \xi, \theta ] \in \theta, \forall \eta \in \Theta \}$ of the {\it normalizer} of $\Theta$ in $T$ is useless in actual practice and must be replaced by the formal definition of ${\tilde{R}}_{q+1}$ given previously through the formal Lie derivative. As such a definition crucially depends on the Spencer operator, it is still not acknowledged. Moreover, as $d {\eta}_{q+1} \in T^* \otimes R_q$, the normalizer can thus be obtained in a purely algebraic way by the condition $ \{ {\tilde{\xi}}_{q+1} , {\eta}_{q+1} \} \in R_q $. In particular, when $n=4$, the Poincar\'e group is of codimension $1$ in its normalizer which is the Weyl group.  \\

\noindent
{\bf EXAMPLE 4.26}: Following Vessiot exactly as we did in the Introduction with $n=1, m=2, k= \mathbb{Q}$, let us consider the generic second order OD equation $y_{xx} - {\omega}^1 y_x + {\omega}^2 y=0$ and copy it twice with $k=1,2$ in order to obtain the linear second order automoprhic system $y^k_{xx} - {\omega}^1 y^k_x + {\omega}^2 y^k=0$ for the standard action of $GL(2)$ on $(y^1,y^2)$ when $k=1,2$. Such a system admits the well known generating Lie form as a quotient of determinants:  \\
\[   {\Phi}^1 \equiv \frac{ 
\mid \begin{array}{ll}
y^1 & y^1_{xx} \\
y^2 & y^2_{xx}
\end{array} \mid }{ 
\mid \begin{array}{ll}
y^1 & y^1_x \\
y^2  & y^2_x
\end{array} \mid } = {\omega}^1 , \hspace{1cm}  
 {\Phi}^2 \equiv \frac{ 
\mid \begin{array}{ll}
y^1_x & y^1_{xx} \\
y^2_x & y^2_{xx}
\end{array} \mid }{ 
\mid \begin{array}{ll}
y^1 & y^1_x \\
y^2  & y^2_x
\end{array} \mid }  ={\omega}^2  \]
In matrix form, the prolongation of the action up to order $q=2$ is:  
\[  ( \begin{array}{lll}
{\bar{y}}^1 & {\bar{y}}^1_x & {\bar{y}}^1_{xx}  \\
{\bar{y}}^2 & {\bar{y}}^2_x & {\bar{y}}^2_{xx}
\end{array}  ) =
( \begin{array}{cc}
a & b \\
c & d 
\end{array} )
( \begin{array}{lll}
y^1 & y^1_x & y^1_{xx}  \\
y^2 & y^2_x & y^2_{xx}                                                                               
\end{array} )  \]
We have the automorphic extension $k \subset K \subset L$ with $K=< {\Phi}^1, {\Phi}^2 > $ and $L = k < y^1,y^2 >$. \\
We obtain therefore $A= {\bar{M}} {M}^{-1}$ in the matrix form:  \\
\[    \bar{M}=A M \Rightarrow  ( \begin{array}{ll}
a & b \\
c & d 
\end{array} )=
( \begin{array}{ll}
{\bar{y}}^1 & {\bar{y}}^1_x \\
{\bar{y}}^2 & {\bar{y}}^2_x
\end{array} )
{( \begin{array}{ll}
y^1 & y^1_x \\
y^2 & y^2_x
\end{array} ) }^{-1} \]
on the condition to have the non-zero Wronskian determinant condition $\Psi \equiv y^1 y^2_x - y^2 y^1_x \neq 0 $.\\
The $4$ infinitesimal generators of the target action are:   \\
\[  \fbox{  $  \Theta= \{ {\theta}_1= y^1 \frac{\partial}{\partial y^1}, {\theta}_2= y^2 \frac{\partial}{\partial y^1}, 
{\theta}_3= y^1 \frac{\partial}{\partial y^2}, {\theta}_4 = y^2 \frac{\partial}{\partial y^2} \}= \{ y^l\frac{\partial}{\partial y^k} \}   $  }   \] 
that can be prolonged at any order like 
${\theta}_1= y^1\frac{\partial}{\partial y^1} + y^1_x \frac{\partial}{\partial y^1_x} + y^1_{xx} \frac{\partial}{\partial y^1_{xx}} + ...$ and so on. \\
One obtains easily the reciprocal distribution at order $1$:  \\
\[ \fbox{  $  \Delta = \{  {\delta}_1= y^k \frac{\partial}{\partial  y^k}, {\delta}_2= y^k_x \frac{\partial}{\partial y^k}, 
{\delta}_3= y^k \frac{\partial}{\partial y^k_x}, {\delta}_4= y^k_x \frac{\partial}{ \partial y^k_x}  \}  $  }  \] 
in such a way that the rank with respect to $(y^1, y^2)$ and $(y^1_x, y^2_x)$ is indeed maximum equal to $2$ provided that the Wronskian determinant does not vanish. Accordingly, one can extend each $\delta$ to $L {\otimes }_K L$ by setting for example ${\delta}_1= y^k \frac{\partial}{\partial y^k} + {\bar{y}}^k\frac{\partial}{\partial {\bar{y}}^k} , ..., {\delta}_4= y^k_x \frac{\partial}{\partial y^k_x}+ {\bar{y}}^k_x \frac{\partial}{\partial {\bar{y}}^k_x}$ and we have:  \\
\[    \Psi=det(M)\neq 0, \,\,\,  \delta {\bar{M}}= (\delta A) M + A (\delta M)= A( \delta M)  \Rightarrow 
\delta A=0  \]
In particular $\fbox{ $a= \frac{y^2_x}{\Psi}{\bar{y}}^1 - \frac{y^2}{\Psi} {\bar{y}}^1_x \in cst(L{\otimes}_K L) $ } $ because ${\delta}_1 \Psi= \Psi, {\delta}_2 \Psi=0, {\delta}_3 \Psi =0, {\delta}_4\Psi= \Psi $ for $SL(2) \lhd GL(2)$ when $ad-bc=1$, a result not evident at first sight. We finally  notice that $A = (\frac{\partial \bar{y}}  { \partial y} )$ and $ \frac{{\partial}^2 \bar{y}}{ \partial y \partial y}=0$ for the differential automorphic extension $L/K$ and we must add $ det(A)=1$ for the differential automorphic extension $L/K'$ with now $K' = K< \Psi>$. We shall explain in the next section why we must take out ${\delta}_1$ and ${\delta}_2$ in the differential algebraic framework if we consider a Lie group as a Lie pseudogroup, that is {\it if we do not inroduce any longer the parameters} $(a,b,c,d)$. We let the reader spend a few minutes now in order to imagine how to manage with such a target.  \\  \\

\noindent
{\bf 5) DIFFERENTIAL ALGEBRA}  \\

The purpose of this section is to revisit the theory of algebraic pseudogroups by using Hopf rings in a way similar to the one pioneered by Bialynicki-Birula in ([2, 3]) (See also [28] for more details and examples). For simplicity, we shall restrict our study to the general situation as the special situation can be treated by restricting the various distributions. For example, the Picard-Vessiot can be achieved for the single OD equation $y_{xx}=0$ by transforming it into the OD automorphic system 
$y^1_{xx}=0, y^2_{xx}=0$. In this case, we may choose $k=K=\mathbb{Q}$ and $L= \mathbb{Q}(y^1,y^2, y^1_x,y^2_x)$ considered as a differential field by setting $d_x y^k = y^k_x, d_x y^k_x=0$ for $k=1,2$ as in the Introduction.. \\

The first comment is to notice that $cst(L)=k(y)$ when $L=k < y >$ and to exhibit a link between between algebraic pseudogroups and the constants of $L {\otimes}_K L$ for $\Delta$ like in the last examples. As no reference can be quoted, we provide a motivating example.   \\

\noindent
{\bf EXAMPLE 5.1}: With $ n=1, m=2, q=1, k=\mathbb{Q}$ and $\Phi\equiv y^2 y^1_x$, let us consider the differential automorphic extension $L/K$ with $K=k< \Phi > \subset  L =k< y^1, y^2 >$. The corresponding algebric pseudogroup is defined by $ {\bar{y}}^2 {\bar{y}}^1_x=y^2 y^1_x \Rightarrow d {\bar{y}}^1 \wedge d {\bar{y}}^2 = d y^1 \wedge d y^2$ or, equivalently, by the first order involutive system of finite Lie equations:  \\
\[  \frac{\partial{\bar{y}}^1}{\partial y^1}= \frac{y^2}{{\bar{y}}^2}= \frac{ {\bar{y}}^1_x}{y^1_x}, 
    \frac{\partial {\bar{y}}^2}{\partial y^2}= \frac{{\bar{y}}^2}{y^2} = \frac{y^1_x}{{\bar{y}}^1_x},
    {\bar{y}}^2_x= \frac{\partial {\bar{y}}^2}{\partial y^1} y^1_x + \frac{\partial {\bar{y}}^2}{\partial y^2} y^2_x \Rightarrow 
    \frac{\partial {\bar{y}}^2}{\partial y^1}= \frac{{\bar{y}}^2_x}{y^1_x} - \frac{y^2_x}{{\bar{y}}^1_x} \]
Starting with ${\delta}_1= y^k_x \frac{\partial }{ \partial y^k_x} $ and ${\delta}_2=y^2 \frac{\partial }{\partial y^2_x}$ and extending these derivations from $L$ to $ L {\otimes}_K L$, we discover that ${\delta}_1$ and ${\delta}_2$ kill the $ 4$  element of the matrix $\frac{\partial {\bar{y}}}{\partial y}$ though this fact is not evident for:  \\
\[  {\delta}_2 (\frac{\partial {\bar{y}}^2}{\partial y^1})= \frac{{\bar{y}}^2}{y^1_x} - \frac{y^2}{ {\bar{y}}^1_x}=
\frac{{\bar{y}}^2{\bar{y}}^1_x - y^2 y^1_x}{y^1_x {\bar{y}}^1_x}= 0  \]
We obtain therefore $k(y){\otimes}_k k({\bar{y}})\in cst (L {\otimes }_K L)$ and $ \frac{\partial {\bar{y}}}{\partial y} \in cst (L {\otimes}_K L)$. \\
Prolonging to order $2$ is even more tricky with, for example:  \\
\[  {\bar{y}}^1_{xx}= \frac{\partial {\bar{y}}^1}{\partial y^1} y^1_{xx} + \frac{{\partial}^2 {\bar{y}}^1}{\partial y^1 \partial y^1}(y^1_x)^2 \Rightarrow  \frac{{\partial}^2 {\bar{y}}^1}{\partial y^1 \partial y^1} =
\frac{1}{(y^1_x)^2}{\bar{y}}^1_{xx} - \frac{y^2y^1_{xx}}{(y^1_x)^2} \frac{1}{{\bar{y}}^2}  \]
and we let the reader check that this term is killed by both ${\delta}_1, ..., {\delta}_4$. Similarly, we get:  \\
\[  {\bar{y}}^2_{xx}= \frac{\partial {\bar{y}}^2}{\partial y^1} y^1_{xx} + \frac{{\bar{y}}^2}{y^2}y^2_{xx} +
\frac{\partial {\bar{y}}^2}{\partial y^1 \partial y^1}(y^1_x)^2 + \frac{1}{y^2} \frac{\partial {\bar{y}}^2}{\partial y^1} y^1_x y^2_x + d_x(\frac{{\bar{y}}^2}{y^2})y^2_x \]
Applyig ${\delta}_3= y^2 \frac{\partial}{\partial y^2_x} + 2 y^2_x \frac{\partial }{ \partial y^2_{xx}}$ while taking into account the previous result at order $1$, we get:  \\
\[  2 {\bar{y}}^2_x = {\delta}_3 (\frac{{\partial}^2 {\bar{y}}^2}{\partial y^1 \partial y^1})(y^1_x)^2 + \frac{\partial {\bar{y}}^2}{\partial y^1}y^1_x + \frac{{\bar{y}}^2}{y^2}y^2_x + {\bar{y}}^2_x   \]
Using the relation $ \frac{\partial {\bar{y}}^2}{\partial y^1}y^1_x + \frac{{\bar{y}}^2}{y^2}y^2_x ={\bar{y}}^2_x  $, we finally obtain $ {\delta}_3(\frac{{\partial}^2 {\bar{y}}^2}{\partial y^1 \partial y^1}) = 0 $ 
and so on.  \\

Though the extension of reciprocal distributions can be done by induction on the order (See [28] for details), we do not know any reference on their explicit computation that may be quite difficult as shown by the next example.  \\

\noindent
{\bf EXAMPLE 5.2}: With $m=n=1, k=\mathbb{Q}$ but $q=3$, we have:  \\
\[  \bar{y}_x=\frac{\partial \bar{y}}{\partial y} y_x, \,\,
\bar{y}_{xx}= \frac{\partial \bar{y}}{\partial y}y_{xx}+ \frac{{\partial}^2\bar{y}}{\partial y^2}(y_x)^2, \,\,
\bar{y}_{xxx}= \frac{\partial \bar{y}}{\partial y} y_{xxx} + 3 \frac{{\partial}^2\bar{y}}{\partial y^2}y_x y_{xx} + \frac{{\partial}^3 \bar{y}}{\partial y^3} (y_x)^3  \]
\[ \fbox{ $  \begin{array}{lcl}
\eta & \rightarrow & \frac{\partial}{\partial y}\\
{\eta}_y & \rightarrow & y_x\frac{\partial}{\partial y_x} + y_{xx} \frac{\partial}{\partial y_{xx}} + y_{xxx} \frac{\partial}{\partial y_{xxx}}  \\
{\eta}_{yy} & \rightarrow & (y_x)^2 \frac{\partial}{\partial y_{xx}} + 3 y_x y_{xx} \frac{\partial}{\partial y_{xxx}} \\
{\eta}_{yyy} & \rightarrow & (y_x)^3 \frac{\partial}{\partial y_{xxx}} 
\end{array}   $  } \]
The reciprocal distribution is:   \\
\[ \fbox{  $ \begin{array}{lcr}
- {\xi}_x & \rightarrow  & y_x \frac{\partial}{\partial y_x} + 2 y_{xx} \frac{\partial}{\partial y_{xx}} + 3 y_{xxx} \frac{\partial}{\partial y_{xxx}}  \\
 - {\xi}_{xx} & \rightarrow & y_x \frac{\partial}{\partial y_{xx}} + 3 y_{xx} \frac{\partial}{\partial y_{xxx}}  \\
  - {\xi}_{xxx} & \rightarrow & y_x \frac{\partial}{\partial y_{xxx}} 
 \end{array}  $  }  \]

Collecting these results, we obtain the following crucial theorem bringing the need to revisit DGT as in ([3, 28]):\\

\noindent
{\bf THEOREM 5.3}: The groupoid components $ ( \frac{\partial \bar{y}}{\partial y}, \frac{{\partial}^2 \bar{y}}{\partial y^2}, 
\frac{{\partial}3 \bar{y}}{\partial y^3}, ... ) $ up to any order $q$ are constants for the reciprocal distribution up to order $q$ that can be expressed as rational functions of all the jet components $ (y, y_x, y_{xx}, y_{xxx}, ... )$ and $ (\bar{y}, {\bar{y}}_x, {\bar{y}}_{xx}, {\bar{y}}_{xxx}, ... )$ up to order $q$. We obtain therefore the Hopf ring $k[\Gamma] \subset cst (L {\times }_K L)$ with ring of fractions $k(\Gamma) = cst (Q((L{\otimes}_K L))$ as a direct sum of differential fields for the target derivative $d_y$ though $Q(L{\otimes}_K L)$ is a direct sum of differential fields for the source derivative $d_x$ with an isomorphism 
$Q(L {\otimes}_{k(y)} k[\Gamma]) \simeq Q (L {\otimes}_K L) $.  \\

\noindent
{\it Proof}:

$L/k(y)$ is a {\it regular} extension because $k(y)$ is algebraically closed in $L$.

\hspace*{12cm}    $\Box $ \\

The proof of the next lemmas is left to the reader through computations with local coordinates (See [28], p 404 for details):  \\

\noindent
{\bf LEMMA 5.4}: If the vector field $W(q+1)= {\Sigma}_{1 \leq \mid \nu \mid \leq q+1} a^k_{\nu}(y_{q+1}) \frac{\partial }{\partial y^k_{\nu}}$ are commuting with vertical vector fields $V(q+1)\in \Theta$, we have the relations:  \\
\[ W(q+1) d_i\Phi= d_i (W(q) \Phi) -{\Sigma}_{0 \leq \mid \mu \mid \leq q} (d_ia^k_{\mu} - a^k_{\mu +1_i} ) \frac{\partial}{\partial y^k_{\mu}} \]
\[ [ V(q+1),W(q+1)] d_i\Phi =V(q+1)(W(q+1) d_i\Phi) - W(q+1) (V(q+1) d_i \Phi)= V(q+1)(W(q+1)d_i\Phi)  \]
 \[\Rightarrow W((q+1)d_i \Phi =  L_i(\Phi, d\Phi) \]  
 
 \noindent
 {\bf LEMMA 5.5}: One has the useful formula for prolongations of source transformations:  \\
 \[   \flat ({\xi}_{q+1})d_i\Phi =d_i (\flat({\xi}_q)\Phi) - {\xi}^r_i d_r \Phi  - \flat (d{\xi}_{q+1}({\partial}_i))\Phi  \]

  The following examples explain the origin of the well known {\it Wronskian determinant } that is existing in the classical Picard-Vessiot theory but in a completely different setting (See [28], p 401).   \\
  
  \noindent 
  {\bf EXAMPLE 5.6}:  When $m=n=1$ and $k=\mathbb{Q}$, let us consider the differential automorphic extension $L/K$ with $K=k<\Phi> \subset L=k<y>$ and $\Phi= \frac{y_x}{y} $. The underlying Lie group action is $\bar{y}=ay \Rightarrow {\bar{y}}_x=ay_x$ when $a=cst$ and it is clear that the word "{\it constant}" is not well defined. On the contrary, the underlying Lie pseudogroup is defined over $k$ by the nonlinear first order system $\frac{1}{\bar{y}} \frac{\partial \bar{y}}{\partial y}=\frac{1}{y}$ in Lie form or, equivalently, by the differential algebraic OD equation $y \frac{\partial \bar{y}}{\partial y} - \bar{y}=0$ over the standard target differential field $k(y)$. Differentiating with respect to $y$, we obtain the linear second order equation $\frac{{\partial}^2 \bar{y}}{\partial y^2}=0$ that needs not to be integrated. As for the PHS law, we have at once the two relations $\frac{\partial \bar{y}}{\partial y}= \frac{{\bar{y}}_x}{y_x}= \frac{\bar{y}}{y}, \frac{{\partial}^2 \bar{y}}{\partial y^2}=0$ for second order jets.  \\
  Now, we have the target transformations:  \\
  \[ \theta= y\frac{\partial}{\partial y} \Rightarrow \sharp (j_2(\theta)) ={\rho}_2(\theta) = y\frac{\partial}{\partial y} + y_x\frac{\partial}{\partial y_x} + y_{xx} \frac{\partial}{\partial y_{xx}}  \]
because $j_2(\theta)= (y,1,0)$ over the target.  \\
For the commuting distribution used in the last lemmas, we may use:  \\
\[ \delta (1)=y_x \frac{\partial }{\partial y_x}, \,\,\, \delta (2) = y_x \frac{\partial}{\partial y_x} + 2 y_{xx} \frac{\partial }{\partial y_{xx}} \]
because $ \flat({\xi}_2)= \xi \frac{\partial}{\partial x} - y_x {\xi}_x \frac{\partial}{\partial y_x} - 
(y_x {\xi}_{xx} + 2 y_{xx} {\xi}_x )\frac{\partial}{\partial y_{xx}} $ ({\it care to the factor 2} ) and thus:   \\
\[  - {\xi}_x \rightarrow y_x \frac{\partial}{\partial y_x} + 2 y_{xx} \frac{\partial}{\partial y_{xx}} , \hspace{1cm}   - {\xi}_{xx} \rightarrow y_x \frac{\partial}{\partial y_{xx}}  \]
and we have indeed $ det ( \begin{array}{cc}
y_x & 2y_{xx} \\
0 & y_x
\end{array} ) =(y_x)^2 \neq 0 $.  \\
With $\Phi = \frac{y_x}{y}  \Rightarrow d_x \Phi = \frac{y_{xx}}{y} - (\frac{y_x}{y})^2 = \frac{y_{xx}}{y} - {\Phi}^2 $, we get therefore:  \\
\[ \delta (1) \Phi= \Phi \Rightarrow   \delta (2) d_x\Phi= d_x \Phi + (y_x \frac{\partial \Phi}{\partial y} + y_{xx} \frac{\partial \Phi}{\partial  
y_x} =d_x\Phi + ( - {\Phi}^2 + \frac{y_{xx}}{y})= 2 d_x \Phi \]
We finally notice that $d_x (\frac{\partial \bar{y}}{\partial y})=d_x (\frac{\bar{y}}{y})=\frac{{\bar{y}}_x}{y} - \frac{\bar{y} y_x}{y^2}=
\frac{y_x}{y}(\frac{\partial \bar{y}}{\partial y} - \frac{\bar{y}}{y})=0$.  \\
It follows from the chain rule for derivatives that we have also $d_x(\frac{\partial \bar{y}}{\partial y})=y_x \frac{{\partial}^2\bar{y}}{\partial y^2}=0$ and thus $\frac{{\partial}^2 \bar{y}}{\partial y^2}=0 $ whenever $y_x\neq 0$. This result is coherent with the fact that the OD equation $y_x=0$ is invariant by {\it any} diffeomorphism $\bar{y}=g(y)$ contrary to a linear OD equation of the form $y_x -\omega y=0$ with $\omega\neq 0$ in a differential field.  \\

\noindent
{\bf EXAMPLE  5.7}: When $n=1, m=2$ and $k=\mathbb{Q}$, let us consider the Picard-Vessiot differential automorphic extension $L/K$ with $K=k<{\Phi}^1,{\Phi}^2>\subset L=k<y^1,y^2>$ along Example 4.26. The underlying Lie group action is ${\bar{y}}^1=a y^1 + b y^2, {\bar{y}}^2=c y^1 + d y^2$ or simply $\bar{y}=Ay$ when $(a,b,c,d)$ are constants and it is clear that the word "{\it constant}" is not well defined, like in the previous example. Moreover, the underlying Lie pseudogroup is surely defined over $k$ by the nonlinear first order system $ {\bar{y}}^u= \frac{\partial  {\bar{y}}^u}{\partial y^k} y^k$ or simply $\bar{y}=\frac{\partial {\bar{y}}}{\partial y} y$ over the standard target differential field $k(y)$. However, differentiating these equations with respect to $y^1$ and $y^2$ as in the previous example, we cannot obtain the linear second order equations $\frac{{\partial}^2 y^u}{\partial y^k\partial y^l}=0$ or simply $\frac{{\partial}^2\bar{y}}{\partial y^2}=0$ as before and {\it these second order PD equations must be added independently}. The underlying Lie pseudogroup {\it must} therefore be defined by second order PD equations and we are no longer allowed to exhibit solutions in the classical form $\frac{\partial {\bar{y}}}{\partial y}=A$ with $A=cst$.   \\
Now, we have the $4$ target transformations ${\theta}^l_k=y^l \frac{\partial }{\partial y^k}$ for $k, l=1, 2$ and their $4$ prolongations up to any order. \\
As for the commuting distribution, we have the two ${\delta}_1= y^k\frac{\partial}{\partial y^k_x}, 
{\delta}_2= y^k_x \frac{\partial}{\partial y^k_x}$ at the order one that can be completed by the two ${\delta}_3= y^k\frac{\partial}{\partial y^k_{xx}}, {\delta}_4= y^k_x\frac{\partial}{\partial y^k_{xx}}$ at the order two. Prolonging source transformations as in the previous example, we obtain ({\it care again to the factor 2} ):  \\
\[   - {\xi}_x \rightarrow y^k_x \frac{\partial}{\partial y^k_x} + 2 y^k_{xx} \frac{\partial}{\partial y^k_{xx}}= {\delta}_2  - 2 {\Phi}^2 {\delta}_3  + 2 {\Phi}^1 {\delta}_4 , \hspace{1cm}   - {\xi}_{xx} \rightarrow {\delta}_4=y^k_x \frac{\partial}{\partial y^k_{xx}}  \]
It follows from the chain rule for derivatives that we have also ${\bar{y}}^u_x=d_x(\frac{\partial {\bar{y}}^u}{\partial y^k}y^k)= \frac{{\partial}^2 {\bar{y}}^u}{\partial y^k \partial y^l} y^k y^l_x + \frac{\partial {\bar{y}}^u}{\partial y^k}y^k_x   $ and thus both 
${\bar{y}}^u=\frac{\partial {\bar{y}}^u}{\partial y^k}y^k, \,\,\,        {\bar{y}}^u_x= \frac{\partial {\bar{y}}^u}{\partial y^k}y^k_x$. \\
Extending ${\delta}_1$ to the derivation ${\delta}_1=y^k\frac{\partial}{\partial y^k_x} + {\bar{y}}^k\frac{\partial}{\partial {\bar{y}}^k_x} $  of $L {\otimes}_K L$ while applying it, we get successively:  \\
\[  y^k {\delta}_1 (\frac{\partial {\bar{y}}^u}{\partial y^k})=0, \hspace{1cm}  y^k_x {\delta}_1 (\frac{\partial {\bar{y}}^u}{\partial y^k})=0 \]
The determinant of this $2 \times 2$ linear system for each $u$ is the {\it Wronskian determinant}: \\
\[   det ( \begin{array}{cc}
y^1 & y^2  \\
y^1_x & y^2_x. 
\end{array} ) = y^1 y^2_x - y^2 y^1_x \neq 0 \]  \\
We may similarly extend ${\delta}_2$ to the derivation ${\delta}_2= y^k_x \frac{\partial}{\partial y^k_x} + {\bar{ y}}^k_x \frac{\partial}{\partial {\bar{y}}^k_x}$ of $L {\otimes}_K L$ in order to check that ${\delta}_2 $ also kills $ \frac{\partial {\bar{y}}}{\partial y}$ whenever $y^1 y^2_x - y^2 y^1_x \neq 0$.   \\

The extension of these results to an arbitrary $m$ is elementary and left to the reader.   \\

\noindent
{\bf EXAMPLE 5.8}: With $k=\mathbb{Q}, m=2,n=1$ let us consider the general automorphic exension $k \subset K \subset K_0 \subset L$ for the Lie groupoid of isometries of the Euclidean metric for ${\mathbb{R}}^2$. In this case, taking the determinant of the isometry, we obtain ${\Delta}^2=1$ with $\Delta=\frac{\partial({\bar{y}}^1,{\bar{y}}^2)}{\partial (y^1,y^2)}$, that is $\Delta=\pm 1$. The prolongations at order two of the infinitesimal target transformations are:  \\
\[  \fbox{  $  {\theta}_1= \frac{\partial}{\partial y^1}, {\theta}_2= \frac{\partial}{\partial y^2}, {\theta}_3=  y^1 \frac{\partial}{\partial y^2} - y^2 \frac{\partial}{\partial y^1} + 
y^1_x \frac{\partial}{\partial y^2_x} - y^2_x \frac{\partial}{\partial y^1_x} + y^1_{xx} \frac{\partial}{\partial y^2_{xx}} - y^2_{xx} \frac{\partial}{\partial y^1_{xx}}  $  }  \]
Accordingly, the only generating differental invariant of order one is $\Omega= ({y^1_x})^2 + ({y^2_x})^2$ while the generating differential invariants at strict order two are 
$\Gamma= \frac{1}{2}d_x\Omega= y^1_xy^1_{xx}+ y^2_xy^2_{xx}$ and $\Upsilon=(y^1_{xx})^2 + (y^2_{xx})^2$. We obtain therefore the general differential automorphic extension $K\subset L$ with $K=k<\Omega, \Gamma, \Upsilon>$ and $L=k<y^1,y^2>$. Setting now $\Sigma = y^1_xy^2_{xx} - y^2_xy^1_{xx}$, we may introduce the other intermediate differential field $K_0= k< \Omega, \Gamma, \Sigma>$ with $dim_K(K_0)=\mid K_0/K\mid =2$. \\
Indeed, using jet coordinates, we have: \\
\[   \left( \begin{array}{cc}
{\bar{y}}^1_x & {\bar{y}}^1_ {xx} \\
{\bar{y}}^2_x & {\bar{y}}^2_{xx}  
\end{array} \right) = \left( \begin{array}{cc}
a & b \\
c & d 
\end{array} \right)  \left(  \begin{array}{c c }
y^1_x &y^1_{xx} \\
y^2_x & y^2_{xx} 
\end{array}  \right)    \]
Hence, taking the determinants, we finally obtain $ \bar{\Sigma}= (ad-bc)\Sigma$ with Jacobian $\Delta= (ad-bc)$. Meanwhile, we have the algebraic relation ${\Sigma}^2 + {\Gamma}^2 - \Omega \Upsilon = 0$ and $K_0$ is the algebraic closure of $K$ in $L$. In the present situation, we have an action of the Lie group $G$ on $(y^1,y^2) $ with $ {\bar{y}}^1=ay^1 + by^2 +\alpha, {\bar{y}}^2= cy^1 + d y^2 + \beta$ where $(a,b,c,d,\alpha,\beta)=(1\,\, rotation + 2 \,\, translations)$ are ordinary constants. The connected component $G^0$ of the identity determines the differential automorphic extension $L/K_0$. In the case of $G$, the symbol $g_2$ of order two is determined by the two linear equations:  \\
\[y^1_x v^1_{xx} +y^2_xv^2_{xx}=0, \,\,\, y^1_{xx}v^1_{xx} + y^2_{xx}v^2_{xx}=0 \]
and $g_2=0$ if and only if $\Sigma \neq 0$. In the case of $G^0$, we have the two linear equations:  \\
\[  y^1_x v^1_{xx} +y^2_xv^2_{xx}=0, \,\,\, y^2_xv^1_{xx} - y^1_x v^2_{xx}=0 \]
and $g_2=0$ if and only if $(y^1_x)^2 + (y^2_x)^2 \neq 0$, in a coherent way with the linearization of the previous algebraic equation. The reciprocal distribution at order two is generated by:  \\
\[   \fbox{  $   {\delta}_1= y^k_x\frac{\partial}{\partial y^k_x} , \,\, {\delta}_2= y^k_{xx}\frac{\partial}{\partial y^k_x}, \,\, {\delta}_3= y^k_x\frac{\partial}{\partial y^k_{xx}}, \,\, 
      {\delta}_4= y^k_{xx}\frac{\partial}{\partial y^k_{xx}}   $  } \]
and is easily seen to stabilize both $K$ and even $K_0$ because we have:  \\
\[      {\delta}_1 \Sigma= \Sigma, \,\,\,{\delta}_2 \Sigma= 0, \,\,\, {\delta}_3 \Sigma = 0, \,\,\,{\delta}_4 \Sigma = \Sigma   \]     
and is full rank, equal to $4$, whenever $\Sigma \neq 0$. However, $K_0$ is algebraically closed in $L$ and thus $k[G^0]\subset cst(L{\otimes}_{K_0} L)$ is an integral domain because $L{\otimes}_{K_0} L$ is an integral domain. On the contrary, $L{\otimes}_K L$ is the direct sum of two integral domains because $K_0{\otimes}_K {K_0 }\simeq K_0 \oplus K_0\simeq K_0 {\otimes}_{\mathbb{Q}}(\mathbb{Q} \oplus \mathbb{Q})$ is a direct sum of two fields as $K_0/K$ is even a classical Galois extension with Galois group defined by ${\Delta}^2=1$ and $L/K_0$ is a regular extension, that is $K_0$ is algebraically closed in $L$. With more details, we have $\bar{\Sigma}=\Delta \Sigma \Rightarrow {\bar{\Sigma}}^2= {\Sigma}^2$ that can also be obtained by substraction from the previous algebraic identity and thus ${\bar{\Sigma}}^2 - {\Sigma}^2= (\bar{\Sigma} - \Sigma)(\bar{\Sigma} + \Sigma)=0$ as a way to split the tensor product when ${\Sigma}^2 \in K$ ([28], Remark 4.57 p 115,[32], p 268-271,[36], p 122). Extending each reciprocal distribution $\delta$ of $L/K$ to $ L{\otimes }_K L $ as usual, we let the reader prove directly that $\Delta=\bar{\Sigma}/\ \Sigma \Rightarrow \delta \Delta =0$. We let also the reader repeat the revious computations in the special case by introducing the nonlinear equations $\Omega=\omega \in K, \Gamma = \gamma \in K, \Upsilon = \upsilon \in K, \Sigma = \sigma \in K_0 $ with CC $ \gamma =\frac{1}{2}{\partial}_x \omega, {\sigma}^2= \omega \upsilon - {\gamma}^2 \in K$ whenever $K$ is a given differential field of characteristic zero and $L=Q(K{y}/\mathfrak{p})$ where $\mathfrak{p}$ is a prime differential ideal defining a differential automorphic system along the ideas of Vessiot.  \\  \\

\noindent
{\bf REMARK 5.9}: It is quite difficult to discover the confusion done by Drach, Kolchin and followers between "{\it maximum ideals}" and "{\it prime ideals} when dealing with the 
so-called concept of {\it reductibility} in the Picard-Vessiot or Drach-Vessiot theories. As it is not essential to understand it because it is covered, in any case, by the confusion between "{\it differential algebraic groups}" and "{\it algebraic Lie pseudogroups}" as explained in the Introduction, we only indicate for the interested reader how to recover it by reading backwards the book of Kolchin ([17]) as follows: \\
VI.5: Proposition 13, p 412 $\Rightarrow$ IV; 5: Corollary 2 to Proposition 2, p 152 $\Rightarrow$ \\III.10: Propositions 6 and 7, p 142 $\Rightarrow$ II.1: Theorem 1, p 86.  \\  
This is a hard task indeed in view of the technical content of these results and their proofs. Also, the assumption made by Kolchin that the field $C$ of differential constants, which must be the same for $K$ {\it and} $L$, must be algebraically closed is in contradiction with the spirit of the Galois theory. Another way to gasp at the importance of maximal ideals is to look at the reference([20]) while comparing it to ([3]). It is also interesting to compare this paper to other Hopf / Galois tentatives like in ([8, 13, 25]).      \\  \\

\noindent
{\bf 6) CONCLUSION   } \\

Life is sometimes a funny comedy on the condition to have a good guardian angel !.  \\
While the author of this paper was correcting the proofs of his 1983 GB book "{\it Differential Galois Theory} ", he had to look for one volume of the Japanese Encyclopedic Dictionary of Mathematics when correcting the pages dealing with the nonlinear Spencer operator.  The two thick volumes being on a shelf in his home, a nearby book called "{\it Th\'{e}orie des Corps D\'{e}formables }" and written by the brothers E. and F. Cosserat in 1909 ([9, 31, 33, 37]) fell down, staying opened at the page 137 in which the two brothers were describing the so-called "{\it Cosserat Couple-Stress Equations} " for continuum mechanics as a better substitute for the " {\it Cauchy Stress Equations} ". The author suddenly understood that the Cosserat brothers had only been computing the nonlinear Spencer operator for the group of rigid motions (3 translations + 3 rotations) acting on 3-dimensional space both with its formal adjoint ([secret]). In a few seconds, the comparison provided the third GB book "{\it Lie Pseudogroups and Mechanics} "  published in 1988 ([29]). This chance has also given him the motivation for applying Lie pseudogroups to mathematical physics ([37]) and this comment, contrary to what could be imagined, is {\it surely not out of the scope of this paper}, as it brings heavy doubts on the usefulness of introducing the "{\it Differential Groups} " of Ritt, Kolchin and all followers, despite the early advertising of the author as we already said in the Introduction.  \\   \\

\newpage

\noindent
{\bf 7) REFERENCES}  \\

\noindent
[1] Artin, E.: Galois Theory, Notre Dame Mathematical Lectures, 2, 1942, 1997.  \\
\noindent
[2] Bialynicki-Birula, A.: On the field of rational functions of Algebraic Groups, Pacific Journal of Mathematics, 11 (1961) 1205-1210.  \\
\noindent
[3] Bialynicki-Birula, A.: On Galois Theory of Fields with Operators, Amer. J. Math., 84 (1962) 89-109.  \\
\noindent
[4] Cartan, E.- Einstein, A.: Letters on Absolute Parallelism, 1929-1932, Princeton University Press, Acad\'{e}mie Royale de Belgique, 1979.  \\
\noindent
[5] Cassidy, P. J.: Differential Algebraic Groups, American Journal of Mathematics, 94, 3 (1972) 891-954 (https://doi.org/10.2307/2373764 ) . \\
\noindent
[6] Chase, S. and Sweedler, M.: Hopf Algebras and Galois Theory, Lecture Notes in Mathematics 97, Springer, 1966.  \\
\noindent
[7] Conrad, K.: 1) Applications of Galois theory 2) Galois groups as permutation groups 3) Galois correspondence theorems 4) Galois groups of cubics and quartics 
5) Recognizing Galois groups $S_n$ and $A_n$: Expositary Papers ( https://kconrad.math.uconn.edu ). \\
\noindent
[8] Cresto, T.; Rio, A.; Vela M.: From Galois to Hopf Galois: Theory and Practice, \\
(arXiv: 1403.6300). \\
\noindent
[9] Cosserat, E. and F.: Th\'{e}orie des Corps D\'{e}formables, Hermann, Paris, 1909.  \\
\noindent
[10] Drach, J.: Th\`{e}se de Doctorat, Ann. Ec. Normale Sup. (3) 15 (1898) 243-384.   \\
\noindent
[11] Dummit, D.S. and Foote, R.M.: Abstract Algebra, Selected Exercises by Brian Felix, Third edition, Wiley, 2011, 2017.  \\
\noindent
[12] Freese, R.: 3 examples involving small Galois Groups, Math 611-612, 2014-2015, \\
(https://math.Hawaii.edu/$\sim$ralph ).  \\
\noindent
[13] Greither, C. and Pareigis, B.: Hopf Galois Theory for separable Field extensions, Journal of Algebra, 106, 1 (1987) 239-258.  \\
\noindent
[14] Jacobson, N.: Lectures in Abstract Algebra, Vol VI, Theory of Fields and Galois Theory, Van Nostrand, 1964.  \\
\noindent
[15] Janet, M.: Sur les Syst\`{e}mes aux D\'{e}riv\'{e}es Partielles, Journal de Math., 8 (1920) 65-121.  \\
\noindent
[16] Kaplanski, I.: An Introduction to Differential Algebra, Hermann, 1957, 1976.  \\
\noindent
[17] Kolchin, E.R.: Differential Algebra and Algebraic groups, Academic Press, 1973.  \\
\noindent
[18] Kolchin, E.R.: Differential Algebraic Groups, Pure and Applied Mathematics, 114, Academic Press, 1985.  \\
\noindent
[19] Kovacic, Geometric Characterization of Strongly Normal Extensions, Transactions of the American Mathematical Society, 9 (2006) 4136- 4157).  \\                                                                                                                                                                                                                                                                                                                                                                                                                                                                                                                 
\noindent                                                                                                                                                                   
[20] Kovacic, J. J.: Picard-Vessiot Theory, Algebraic Groups and Group Schemes, 2005) (Kolchin Seminar in Differential Algebra).  \\ 
\noindent
[21] Kumpera  A, Spencer DC, Lie Equations, {\it Ann. Math. Studies } 73, Princeton University Press, Princeton, 1972.  \\
\noindent
[22] Mac Lane, S. and Birkhoff, G.: Algebra, AMS Chelsea Publishing, Third edition, Solutions D\'{e}velopp\'{e}es des Exercices, A. Mezard, CH. Delorme, Ch. Lavit, J.-C. Raoult, Gauthiers-Villars, Bordas, Paris, 1976. \\
\noindent
[23] Manasse-Boetanni, J.: An Introduction to Galois Theory, 2013, \\
(https://math.uchicago.edu/manasse-boetani) . \\
\noindent
[24] Milne, J.S.: Fields and Galois Theory, Version 5.10, 2022 (www.jmilne.org/math/1996-2022.  \\
\noindent
[25] Montgomery, S.: Hopf Galois Theory: A survey, Geometry and Topology Monographs, 16 (2009) 367-400, 
(DOI: 10.2140/gtm.2009.16.367). \\
\noindent
[26] Nichols, W. and Weisfeiler, B.: Differential Formal Groups of J. F. Ritt, American Journal of Mathematics, 104, 5 (1982) 943-1003 (https://doi.org/10.2307/2374080). \\  
\noindent
[27] Pommaret, J.-F.: Systems of Partial Differential Equations and Lie Pseudogroups, Gordon and Breach, New York, 1978 
(Russian translation by MIR, Moscow, 1983). \\
\noindent
[28] Pommaret, J.-F.: Differential Galois Theory, Gordon and Breach, New York, 1983.\\
\noindent
[29] Pommaret, J.-F.: Lie Pseudogroups and Mechanics, Gordon and Breach, New York, 1988.  \\
\noindent
[30] Pommaret, J.-F.: Partial Differential Equations and Group Theory: New Perspectives for Applications, Kluwer, 1994 (http://dx.doi.org/10.1007/978-94-017-2539-2 ). \\
\noindent
[31] Pommaret, J.-F.: Fran\c{c}ois Cosserat and the Secret of the Mathematical Theory of Elasticity, Annales des Ponts et chauss\'{e}es, 82 ( 1997) 59-66 (Translation by D. H. Delphenich). \\
\noindent 
[32] Pommaret, J.-F.: Partial Differential Control Theory, Kluwer, 2001 (957 pp).\\
\noindent
[33] Pommaret, J.-F.: Parametrization of Cosserat Equations, Acta Mecchanica, 215 (2010) 43-55. \\
(https://doi.org/10.1007/s00707-010-0292-y ) \\
\noindent
[34] J.-F. POMMARET: Spencer Operator and Applications: From Continuum Mechanics to Mathematical Physics, in "Continuum Mechanics-Progress in Fundamentals and Engineering Applications", Dr. Yong Gan (Ed.), ISBN: 978-953-51-0447--6, InTech, 2012, Available from: \\
http://www.intechopen.com/books/continuum-mechanics-progress-in-fundamentals-and-engineerin-applications/spencer-operator-and-applications-from-continuum-mechanics-to-mathematical-physics  \\
\noindent
[35] Pommaret, J.-F.: Deformation Theory of Algebraic and Geometric Structures, Lambert Academic Publisher (LAP), Saarbrucken, Germany, 2016.  \\
\noindent
[36] Pommaret, J.-F.: New Mathematical Methods for Physics, Mathematical Physics Books, NOVA Science Publisher, New York, 2018.  \\
\noindent
[37] Pommaret, J.-F.: Minimum Parametrization of the Cauchy Stress Operator, , Journal of modern Physics, 12 (2021) 453-482, . (https://doi.org/10.4236/jmp.2021.12131106).  \\
\noindent
[38] Pommaret, J.-F.: How Many Structure Constants do Exist in Riemannian Geometry ?, Mathematics in Computer Science, 16, 23 (2022) 
(https://doi.org/10.1007.s11786-022-00546-3).  \\
\noindent
[39] Ritt, J.F.: Differential Algebra, Dover, 1950, 1966.  \\
\noindent
[40] Rotman, J. J.: An Introduction to Homological Algebra, Pure and Applied Mathematics, Academic Press, 1979.  \\
\noindent
[41] Serre, J.-P.: Topics in Galois Theory, Course at Harvard University, Fall 1988, Notes written by Henry Darmon, 1991.  \\
\noindent
[42] Spencer, D. C.: Overdetermined Systems of Partial Differential Equations, Bull. Am. Math. Soc., 75 (1965) 1-114.  \\
\noindent
[43] Stewart, I.: Galois Theory, Chapman and Hall, 1973.  \\
\noindent
[44] Sweedler, M.: Hopf Algebras, Benjamin, New York, 1969.  \\
\noindent
[45] Vessiot E. Sur la Th\'{e}orie des Groupes Infinis, {\it Ann. Ecole Normale Sup.} (1903) 20: 411- 451. (http://numdam.org/ ).  \\
\noindent
[46] Vessiot, E.: Sur la Th\'{e}orie de Galois et ses Diverses G\'{e}n\'{e}ralisations, Ann.Ec. Normale Sup., 21 (1904) 9-85 (http://numdam.org/ ).  \\
\noindent
[47] Zariski, O., Samuel, P.: Commutative Algebra, Van Nostrand, 1958.  \\

\end{document}